\documentclass[a4paper,12pt]{amsart}
\usepackage{amssymb,amscd}
\usepackage{calrsfs}
\usepackage{bm}
\theoremstyle{plain}
\newtheorem{thm}{Theorem}[section]
\newtheorem{prop}[thm]{Proposition}
\newtheorem{lemm}[thm]{Lemma}
\newtheorem{coro}[thm]{Corollary}

\theoremstyle{definition}

\theoremstyle{remark}
\newtheorem{rem}[thm]{Remark}
\newtheorem*{note}{Note}

\newcommand{\U}{\mathcal U}
\newcommand{\T}{\mathcal T}
\newcommand{\W}{\mathcal W}

\newcommand{\field}[1]{\mathbb{#1}}
\newcommand{\C}{\field{C}}

\newcommand{\Q}{\field{Q}}
\newcommand{\R}{\field{R}}
\newcommand{\Z}{\field{Z}}

\DeclareMathOperator{\rank}{rank} \DeclareMathOperator{\tr}{tr}

\DeclareMathOperator{\ch}{ch}
\DeclareMathOperator{\ind}{ind}

\newcommand{\lam}{\lambda}
\newcommand{\var}{\varphi}

\newcommand{\hatvar}{\hat{\varphi}}

\newcommand{\hatV}{\hat{V}}
\newcommand{\hatU}{\hat{U}}

\newcommand{\hatF}{\hat{F}}
\newcommand{\hatH}{\hat{H}}
\newcommand{\hatN}{\hat{N}}
\newcommand{\hatX}{\hat{X}}
\newcommand{\hatW}{\hat{W}}
\newcommand{\hatT}{\hat{T}}

\newcommand{\hatp}{\hat{p}}
\newcommand{\hatx}{\hat{x}}
\newcommand{\hatLam}{\hat{\Lambda}}
\newcommand{\hatlam}{\hat{\lambda}}
\newcommand{\hatpi}{\hat{\pi}}
\newcommand{\HatLam}{\Hat{\Hat{\Lambda}}}
\newcommand{\Hatlam}{\Hat{\Hat{\lambda}}}
\newcommand{\HatF}{\Hat{\Hat{F}}}
\newcommand{\HatN}{\Hat{\Hat{N}}}

\newcommand{\hatgam}{\hat{\gamma}}
\newcommand{\calC}{\mathcal{C}}
\newcommand{\hatcalC}{\hat{\mathcal{C}}}

\newcommand{\brvar}{\breve{\varphi}}
\newcommand{\brf}{\breve{f}}
\newcommand{\brT}{\breve{T}}

\def\l{\langle}
\def\r{\rangle}

\newcommand{\img}{\sqrt{-1}}

\title{Orbifold elliptic genera and rigidity}
\date{}
\author{Akio Hattori}
\address{Graduate School of Mathematical Science, University of Tokyo,
Tokyo, Japan}
\email{hattori@ms.u-tokyo.ac.jp}
 
\begin{document}
\maketitle

\section{Introduction}
\label{sec:intro}

Let $M$ be a closed almost complex manifold on which a compact connected 
Lie group $G$ acts non-trivially. If the first Chern class of $M$ is 
divisible by an integer $N$ greater than 1, then its equivariant 
elliptic genus $\var(M)$ of level $N$ is rigid, i.e., it is constant
as a virtual character of $G$. This result was predicted by 
Witten \cite{W} and proved by Taubes \cite{T}, Bott-Taubes \cite{BT} and 
Hirzebruch \cite{Hi}. Elliptic genus can be defined even for 
almost complex orbifolds. Moreover another genus called orbifold
elliptic genus is defined for orbifolds. A natural question is whether 
the rigidity property holds for these genera on orbifolds or not. It turns
out that the answer is no in general. In \cite{HM} we were concerned 
with related topics. 

In this note a modified orbifold elliptic genus of level $N$ 
will be defined for closed almost complex orbifolds such that 
$N$ is relatively prime to the orders of all isotropy groups. 
One of main results, Theorem \ref{thm:brvar}, states that the 
modified orbifold elliptic genus 
$\brvar(X)$ of level $N$ of an almost complex orbifold $X$ of 
dimension $2n$ such that $\Lambda^nTX=L^N$ for some 
\emph{orbifold} line bundle $L$  
is rigid for non-trivial $G$ action. 
As to the orbifold elliptic genus itself Theorem \ref{thm:hatvar} states that 
the orbifold elliptic genus $\hatvar(X)$ of level $N$ of $X$ 
is rigid for non-trivial $G$ action if  
$\Lambda^nTX=L^N$ for some \emph{genuine} $G$ line bundle $L$. 
Furthermore the orbifold elliptic genus $\hatvar(X)$ 
is rigid for non-trivial $G$ action if  
$\Lambda^nTX$ is trivial as an orbifold line bundle 
(Theorem \ref{thm:hatvartorsion}). 
The last result is essentially due to
Dong, Liu and Ma \cite{DLM}.

Liu \cite{L} gave a proof of rigidity by using modular property 
of elliptic genera for manifolds. 
Our proof of the rigidity for the genera $\brvar(X)$ and $\hatvar(X)$ 
also uses Liu's method.

The organization of the paper is as follows. In Section 2 
we review basic materials concerning orbifolds in general. 
The notion of sectors 
is particularly relevant for later use. In Section 3
we give the definitions of orbifold elliptic genus and modified 
orbifold ellitic genus and the main theorems are stated here. 
Section 4 is devoted to exhibiting fixed point formulae for the above 
genera. The proof of the main results will be given in 
Section 5 and Section 6. In Section 6 some additional results 
related to vanishing property are 
given. Main results in this section are Propositions 
\ref{prop:vanishing}, \ref{prop:hatvanishing} and 
\ref{prop:trivialcone}. 
Section 7 concerns the orbifold $T_y$ genus and its modified one. 
They are always rigid for non-trivial actions of campact connected 
Lie groups and take special forms when 
the orbifold elliptic genera vanish. In Section 8 the 
generalization to the case of stably almost complex orbifolds are 
discussed and it will be shown that main results in Section 3, 
Section 6 and Section 7 also hold for stably almost complex orbifolds. 

The author is grateful to M. Furuta, M. Futaki, A. Kato, Y. Mitsumatsu 
and other members of Furuta's Seminar who attended his talks and 
gave him several useful comments. The presentation of this paper 
was largely improved by their help. He is also grateful to 
M. Masuda , coauthor of the related paper \cite{HM}, for his 
collaboration which initiated the work of present paper.

\section{Sectors}

We first recall some basic facts about orbifolds. 
We refer to \cite{S} and \cite{D} for relevant  
definitons and to \cite{HM} for notations used here. 
In \cite{S} orbifolds were called $V$-manifolds.

Let $X$ be a closed orbifold of dimension $n$ and let $\U$ be an atlas of $X$. 
$\U$ is a collection of (orbifold) charts $\{(V_\mu,U_\mu,H_\mu,p_\mu)\}$ 
where $U_\mu$ is an open set in $X$, $V_\mu$ is a smooth manifold 
of dimension $n$, $H_\mu$ is a
finite group acting on $V_\mu$ and $p_\mu$ is a map $V_\mu \to U_\mu$ inducing 
a homeomorphism from $V_\mu/H_\mu$ onto $U_\mu$. The collection $\{U_\mu\}$ 
is assumed to contain a basis of neighborhoods for each point $x\in X$. 
Let $(V_{i},U_{i},H_{i},p_{i}), \ i=1,2,$ be two
orbifold charts in $\U$ such that $U_1\subset U_2$. 
A pair of injective group homomorphism $\rho:H_{1} \to H_{2}$
and $\rho$-equivariant open embedding $\psi:V_{1} \to V_{2}$ 
covering the inclusion map $U_{1}\to U_{2}$ with the 
following property
\begin{equation}\label{eq:1}
 \{h\in H_{2} \mid h(\psi(V_{1}))\cap \psi(V_{1})
  \not=\emptyset\}=\rho(H_{1})
\end{equation}
is called an injection of charts 
and will be written as
\[
 \Psi=(\rho,\psi):(V_{1},U_{1},H_{1},p_{1})\to 
    (V_{2},U_{2},H_{2},p_{2}). \] 
Note that, if $h\in H_{2}$, then 
$(c_h \circ\rho, h\circ\psi)$ is an injetion of charts,
where $c_h:H_{2}\to H_{2}$ denotes the conjugation by $h$.
These charts must satisfy 
the following compatibility condition. 
If $U_{1}$ and $U_{2}$ have a non-empty intersection, then, 
for each $x\in U_{1}\cap U_{2}$, 
there are a chart $(V,U,H,p)\in \U$ 
with $x\in U\subset U_{1}\cap U_{2}$ and 
injections of charts
$ \Psi_i=(\rho_i,\psi_i):(V,U,H,p)\to 
    (V_{i},U_{i},H_{i},p_{i}),\ i=1,2$. 

\begin{note}
We do not assume the action of $H_\mu$ on $V_\mu$ is effective. If 
the effectiveness is assumed then the condition \eqref{eq:1} 
automatically follows. 
\end{note}
When one considers additional structures like Riemannian
structure or almost complex structure, the $V_\mu$ are assumed 
to have the structures in question, and the action of $H_\mu$ and 
the maps $\psi$ are assumed to preserve those structures. 

Let $\Psi=(\rho,\psi):(V_{1},U_{1},H_{1},p_{1})\to 
    (V_{2},U_{2},H_{2},p_{2})$ be an injection of charts. 
One readily sees that the isotropy subgroup $H_{1,v}$ at a point 
$v\in V_{1}$ is isomorphic to that of 
$H_{2}\in V_{2}$ at $\psi(v)$ 
as a consequence of \eqref{eq:1}. Thus the isomorphism class 
of the group $H_{1,v}$ depends only on $x=p_1(v)\in X$. 
It is called the \emph{isotropy group} of the 
orbifold $X$ at $x$ and will be denoted by $H_x$. 
  
For each $x\in X$ there is an orbifold chart $(V_x,U_x,H_x,p_x)$ such that
$p^{-1}(x)$ is a single point, and hence $H_x$ is the isotropy group at $x$. 
Moreover we can take $U_x$ small enough so that 
$V_x$ is modelled on an open ball in $\R^n$ with a linear $H_x$ action. 
Then the fixed point set $V_x^h$ of each $h\in H_x$ is connected. 
Such a chart will be called a \emph{reduced chart} centered at $x$.  
Hereafter we sometimes write simply $(V,U,H)$ for an orbifold 
chart $(V,U,H,p)$, when the meaning of $p$ is clear from the context.
We also assume that orbifolds considered hereafter are closed 
(compact, connected and without boundary) orbifolds.

Fix a finite group $H$ and set $X_H=\{x\in X \mid H_x\cong H\}$. 
Then $X_H$ is a smooth manifold. In fact, if $(V_x,U_x,H_x,p_x)$ 
is a reduced chart centered at $x\in X_H$, then $p_x$ induces a 
homeomorphism between $V_x^{H_x}$ and $X_H\cap U_x$. 
A connected component of $X_H$ is called a \emph{stratum} of $X$. The
totallity of strata is called a stratification of $X$. If $X$ is 
connected then there is a unique stratum $S_H$ such that $H$ is the
minimum with respect to the obvious ordering among isotropy groups 
induced by inclusions. This stratum is 
called principal stratum of $X$ and the order $|H|$ is called the 
multiplicity of the connected orbifold $X$ and is denoted by $m(X)$. 
The principal stratum is open and dense in $X$. 
In case the actions of isotropy groups are effective we have 
$m(X)=1$. 

A map $f:X\to X'$ from an orbifold $X$ to another orbifold $X'$ is 
called \emph{smooth} if, for each orbifold 
chart $(V,U,H,p)$ around $x$ and $(V',U',H',p')$ around $f(x)$ 
such that $f(U)\subset U'$, there is a pair $(\rho,\psi)$ of group 
homomorphism $\rho:H\to H'$ and $\rho$-equivariant smooth 
map $\psi:V\to V'$ satisfying the relation 
$p'\circ \psi=f\circ p$. The totality of these pairs $(\rho,\psi)$ 
is required to satisfy the obvious compatibility relation with respect 
to injections of orbifolds charts for $X$ and $X'$.
Such a pair $(\rho,\psi)$ is called a map of 
orbifold chart covering $f$ and will be written 
$(\rho,\psi): (V,U,H,p)\to (V',U',H',p')$. 

Let $G$ be a Lie group. An action of $G$ on an orbifold $X$ is a
smooth map $G\times X\to X$ which satisfies the usual rule of 
group action.

Let $X'$ be an orbifold with atlas $\U '$. Let $X$ be a subspace 
of $X'$ such that,  
for each chart $(V',U',H,p)\in \U '$ of $X'$, 
$p^{-1}(X\cap U')$ is an $H$-invariant submanifold $V$ of $V'$ 
of dimension $n$. Then $\U=\{(V,X\cap U',H,p|V)\mid (V',U',H,p)\in \U '\}$ 
defines an orbfold structure on $X$. With this structure $X$ is called 
a suborbifold of $X'$.

A triple $(W,X,p)$ of orbifolds $W,X$ and smooth map
$\pi:W\to X$ is called an orbifold vector bundle over $X$ if 
it satisfies the following condition. For any orbifold chart $(V,U,H,p)$ 
for $X$ with $U$ sufficiently small there is an orbifold chart 
$(\tilde{V},\tilde{U},\tilde{H},\tilde{p})$ for $W$ with 
$\tilde{U}=\pi^{-1}(U)$ 
and $\tilde{H}=H$ together with a map of orbifold chart 
$(\rho,\psi):(\tilde{V},\tilde{U},\tilde{H},\tilde{p})\to (V,U,H,p)$ 
covering $\pi:W\to X$ 
such that $\rho=\text{identity}$ and $\tilde{p}:\tilde{V}\to V$ is 
a vector bundle with compatible $H$ action.
The tangent bundle $TX$ of an orbifold $X$ is a typical example of 
vector bundles. Locally its orbifold chart is given by 
$(TV,TV/H,H)$. 
A smooth map $s:X\to W$ is a section if it is locally 
an $H_x$-invariant section. If $s$ is a section then 
$\pi\circ s$ equals the identity map. 
A differential form of degree $q$ is a section of the 
$q$-th exteior power $\Lambda^qT^*X$ of the 
cotangent bundle $T^*X$. 

We are now in a position to give the definition of so called sectors 
of an orbifold $X$. The notion was first introduced by Kawasaki in 
\cite{Ka} in order to describe the index theorem for orbifolds. We set 
\begin{equation*}\label{eq:twistor}
 \hatX =\bigsqcup_{x\in X}\calC(H_x),
\end{equation*}
where $\calC(H)$ denotes the set of conjugacy classes of a group $H$.
We define $\pi: \hatX \to X$ by $\pi(\calC(H_x))=x$. 
An orbifold structure is endowed on $\hatX$ in the following way. 
Let $(V,U,H,p)$ be an orbifold chart for $X$. We set 
\begin{equation*}
 \hatV=\{(v,h)\in V\times H \mid hv=v\}\ \text{and}\ \hatU
  =\bigsqcup_{x\in U}\calC(H_x) \subset \hatX.
\end{equation*}
The group $H$ acts on $\hatV$ by
\[ g(v,h)=(gv,ghg^{-1}). \]
We define $\hatp:\hatV\to \hatU$ by
$\hatp(v,h)=[h]\in \calC(H_{p(v)})$. Then $\hatp$ induces a bijection 
$\hatV/H\to \hatU$. By this bijection we identify 
$\hatV/H$ and $\hatU$. Using this identification one can give a 
topology on $\hatX$ and the collection of quadruples 
$\{(\hatV,\hatU,H,\hatp)\}$ 
defines an orbifold structure on $\hatX$. Moreover, if we define 
$\hat{\pi}:\hatV\to V$ by $\hat{\pi}(v,h)=v$, then 
$(id,\hat{\pi})$ defines a map of orbifold chart 
$(\hatV,\hatU,H,\hatp)\to (V,U,H,p)$. It follows that 
the map $\pi: \hatX \to X$ is a smooth map of orbifolds.

$\hatX$ is not connected unless $X$ is a (connected) 
smooth manifold. Its connected component are sometimes called 
\emph{sectors}. We shall call $\hatX$ the \emph{total sector} of $X$. 
To get reduced charts of $\hatX$ centered at $\hatx\in \hatX$
we proceed as follows. Note that, if $C(h)$ denotes the centralizer 
of $h\in H$, then $C(h)$ acts on the fixed point set $V^h$ of $h$.
If $h$ and $h'$ are conjugate each other in $H$, 
then there is a canonical homeomorphism between $V^h/C(h)$ and 
$V^{h'}/C(h')$. So we can associate the space $V^h/C(h)$ 
to each conjugacy class $[h]$. We also see that there is a 
disjoint sum decomposition of $\hatU=\hatV/H$ 
\[ \hatU=\bigsqcup_{[h]\in \calC(H)}V^h/C(h). \]
Take a reduced orbifold chart $(V_x,U_x,H_x)$ centered at $x$.
Then $V_x^h$ is connected and hence $(V_x^h,V_x^h/C(h),C(h))$ 
gives us a reduce chart of $\hatX$ centered at $[h]\in \calC(H_x)$. 
\begin{note}
The projection $\hatU\to U$ is covered by the inclusion 
map $V_x^h\to V_x$. Hence $\pi: \hatX\to X$ is an immersion of 
orbifold.
\end{note}

When $X$ is connected and its multiplicity $m(X)$ is equal to $1$, 
there is a unique component of $\hatX$ which is mapped
isomorphically onto $X$ by $\pi$, where the point over $x\in X$ is 
the identity element in $H_x$ regarded as an elemnet in $\calC(H_x)$. 
Components with lower dimensions are called \emph{twisted sectors}.

In order to look more closely at a sector, consider a reduced 
chart $(V_x,U_x,H_x,p_x)$ centered at $x$ and take a point $y$ in $U_x$.
$p_x^{-1}(y)$ is an orbit of $H_x$ and the isotropy subgroup $(H_x)_v$ 
of $H_x$ at each point $v\in p_x^{-1}(y)$ is isomorphic to $H_y$. 
This determines an injection of $H_y$ into $H_x$ unique up to 
conjugations by elements of $H_x$. This injection induces a canonical 
map $\rho_{x,y}:\calC(H_y)\to \calC(H_x)$. 

Let $\{\hatX_{\hatlam}\}_{\hatlam\in \hatLam}$ denote the totality 
of sectors. Take $\gamma_x\in \hatX_{\hatlam}$ and take 
a reduced chart $(V_x,U_x,H_x,p_x)$ centered at $x=\pi(\gamma_x)\in X$. 
Then the following lemma gives another local expression of a sector. 
\begin{lemm}\label{lemm:rhoxy}
 \[ \hatX_{\hatlam}\cap \pi^{-1}(U_x)=\{\gamma_y \mid y\in U_x, 
  \gamma_y\in \pi^{-1}(y), \rho_{x,y}(\gamma_y)=\gamma_x \}. \]
\end{lemm}
\begin{proof} 
It is clear from the definition that the left hand side is contained 
in the right hand side. Conversely take an element 
$\gamma_y \in \pi^{-1}(U_x)$ with $\pi(\gamma_y)=y$ and 
$\rho_{x,y}(\gamma_y)=\gamma_x$. Let $v$ be a point in $V_x$ such 
that $p_x(v)=y$, and identify $(H_x)_v$ with $H_y$. If $h\in H_x$ is 
a representative of $\gamma_x$, then $h$ belongs to $(H_x)_v=H_y$ since 
$\rho_{x,y}(\gamma_y)=\gamma_x$. Hence $v\in V_x^h$. This means 
that $\gamma_y=[h]$ belongs to the same component $\hatX_{\hatlam}$ as 
$\gamma_x$. 
\end{proof}
 
We set
\[ \hatH_x =\{h\in H_x \mid V_x^h= V_x^{H_x}\}. \]
This set is closed under conjugations. The conjugacy classes in 
$\hatH_x$ is denoted by $\calC(\hatH_x)$. 
The following equality is immediate from the definition.
\begin{equation*}\label{eq:hatH}
 \calC(\hatH_x)=\{\gamma_x\in \calC(H_x)\mid \gamma_x\not\in 
 \rho_{x,y}(\calC(H_y)) \ \text{for any 
 $y\in U_x\setminus (U_x\cap p_x(V_x^{H_x}))$} \}. 
\end{equation*}

In this paper 
we shall make the assumption \\
(\#) \quad \emph{the fixed point set $V_x^H$ of each subgroup  
$H$ of $H_x$ has even codimension in $V_x$} \\
throughout. 
This is the case of 
stably almost complex orbifolds. 

Let $\hatX_{\hatlam}$ be a sector. 
A point $\gamma_x$ in $\hatX_{\hatlam}$ such that 
$\gamma_x\in \calC(\hatH_x)$ where $x=\pi(\gamma_x)$ 
will be called \emph{generic}. 

\begin{lemm}
The set of all generic points in $\hatX_{\hatlam}$ is a 
connected and dense open set in $\hatX_{\hatlam}$. 
\end{lemm}
\begin{proof}
Let $(V_x^h,V_x/C(h),C(h))$ be a reduced chart centered at
$\gamma_x=[h]$. If $\gamma_x$ lies 
in the principal stratum of $\hatX_{\hatlam}$, then $V_x^h=V_x^{C(h)}$. 
Assume that $V_x^h\not=V_x^{H_x}$ and hence 
$V_x^h\supsetneqq V_x^{H_x}$. Then there exists a poin $x_1\in V_x$ 
such that $C(h)\subset H_{x_1}\subsetneqq H_x$. 

If $V_{x_1}^h\not=V_{x_1}^{H_{x_1}}$ further, then there is a sequence 
of points $x_1,x_2,\ldots \in V_x$ such that 
\[ H_x \supsetneqq H_{x_1}\supsetneqq H_{x_2}\supsetneqq \cdots .\]
Since $H_x$ is a finite group, this sequence terminates at a finite
step and there is a point $y\in V_x$ and $h\in H_y$ such that
$V_y^h=V_y^{H_y}$. Then $\gamma_y=[h]\in \calC(\hatH_y)$ is 
generic. This proves the existence of generic points.

We shall denote by $\hatX_{\hatlam}^{gen}$ the set of generic 
points in $\hatX_{\hatlam}$. If $\gamma_x$ lies in 
$\hatX_{\hatlam}\setminus \hatX_{\hatlam}^{gen}$, then there is a point
$\gamma_{x_0}=[h_0]\in \hatX_{\hatlam}^{gen}$ near $\gamma_x$ such that 
\[ V_x^{h_0}\supsetneqq V_x^{H_x}, \]
where $x=\pi(\gamma_x)$. Since $V_x^{H_x}$ has even codimension in
$V_x^{h_0}$ by the assumption (\#), it follows that 
$V_x^{h_0} \setminus V_x^{H_x}$ is connected. This 
implies that $\hatX_{\hatlam}^{gen}$ is connected. It is also 
open and dense in $\hatX_{\hatlam}$. 
\end{proof}

Let $\{S_\lambda\}_{\lambda\in \Lambda}$ be 
the totality of strata of $X$. 
The isomorphism class of isotropy group at any point in $S_\lambda$ 
depends only on $S_\lambda$ and is denoted by $H_\lambda$. 
The index set $\Lambda$ is a poset 
by the ordering $\prec$ defined by
\[ \lambda\prec \lambda' \Leftrightarrow S_\lambda\subset 
\bar{S}_{\lambda'}. \]
Note that, if $\lambda\prec \lambda'$, then $H_\lambda \supset H_{\lambda'}$
in the sense that there is an injective homomorphism 
$ H_{\lambda'}\to H_\lambda$. 
When $X$ is connected $\Lambda$ has a unique maximal element $\lambda_0$. 
$S_{\lambda_0}$ is the principal stratum and $H_{\lambda_0}$ is the 
minimum isotropy group. 

Since $\pi(\hatX_{\hatlam}^{gen})$ is connected, there is a unique 
$\lambda\in \Lambda$ such that $\pi(\hatX_{\hatlam}^{gen})\subset S_\lambda$. 
This $\lam$ will be denoted by $\pi(\hatlam)$.
Note that $\pi:\hatX\to X$ maps $\hatX_{\hatlam}$ onto $\bar{S}_\lam$. 
The restriction of $\pi$ on 
$\hatX_{\hatlam}^{gen}$ is a covering map onto $S_\lambda$. 

\begin{rem}
An orbifold $X$ is orientable by definition if $V$ can be given an 
orientation and 
the action of $H$ on $V$ is orientation preserving for each orbifold 
chart $(V,U,H)$, and $\psi$ preserves the given orientations of $V$ 
and $V'$ for each
injection $(\rho,\psi):(V,U,H)\to (V',U',H')$. In this case 
$X$ regarded as a $\Q$-homology manifold is orientable. It should be 
noticed however that strata and sectors of $X$ are not 
orientable in general even if the assumption (\#) is satisfied. 
\end{rem}

\section{orbifold elliptic genus}

In this section $X$ will be an almost complex closed orbifold. 
We shall give the definition of elliptic genus $\var(X)$, orbifold 
elliptic genus $\hatvar(X)$ and modified orbifold elliptic genus 
$\brvar(X)$.

If $(V,U,H)$ is an orbifold chart of $X$, then $V$ is an almost complex 
manifold and the action of each element $h\in H$ on $V$ preserves 
almost complex structure. Hence $V^h$ is also an almost complex manifold. 
It follows that the sectors of $X$ are all almost complex orbifolds. 
Let $\hatX_{\hatlam}$ be a sector of $X$, and let $W_{\hatlam}$ be 
the normal bundle of the immersion 
$\pi:\hatX_{\hatlam}\to X$. Take $\gamma_x\in \hatX_{\hatlam}$ 
with $\pi(\gamma_x)=x$. Let $h\in H_x$ be a representative of 
$\gamma_x$. Then $h$ acts on the normal bundle 
$\tilde{W}_x$ of $V_x^h$ in $V_x$. Let 
\begin{equation}\label{eq:eigendecomp}
 \tilde{W}_x=\bigoplus \tilde{W}_{x,i} 
\end{equation}
be the eigen-bundle decompositions with respect to this action of $h$, 
where $h$ acts on $\tilde{W}_{x,i}$ with weight $m_{\hatlam,i}$, 
i.e., by the multiplication of $e^{2\pi\img m_{\hatlam,i}}\not=1$. 
$m_{\hatlam,i}$ is determined modulo integers and depends only on 
$[h]=\gamma_x$. If $y$ is near $x$ 
and $\gamma_y$ lies in $\hatX_{\hatlam}^{gen}$, then 
$\rho_{x,y}(\gamma_y)=\gamma_x$ by Lemma \ref{lemm:rhoxy}. 
Hence $m_{\hatlam,i} \bmod \Z$ is a locally constant function 
of $\gamma_x$ and consequnetly it is constant on $\hatX_{\hatlam}$. 
It follows that the decomposition \eqref{eq:eigendecomp} gives 
the one for the bundle $W_{\hatlam}$:
\begin{equation}\label{eq:eigendecomp2}
 W_{\hatlam}=\bigoplus W_{\hatlam,i}. 
\end{equation}
The $m_{\hatlam,i}$ with $0<m_{\hatlam,i}<1$ is written 
by $f_{\hatlam,i}$. We set 
$f_{\hatlam}=\sum_i f_{\hatlam,i}\dim W_{\hatlam,i}$. 

Let $\tau,\sigma\in \C$ with $\Im(\tau)>0$. We set 
$q=e^{2\pi\img \tau}, \zeta=e^{2\pi\img \sigma}$. 
Let $TX$ be the complex tangent bundle of the almost complex 
orbifold $X$. 
We define formal vector bundles $\T=\T(\sigma), 
\T_{\hatlam}=\T_{\hatlam}(\sigma),
\W_{\hatlam,i}=\W_{\hatlam,i}(\sigma)\ 
 \text{and}\ \W_{\hatlam}=\W_{\hatlam}(\sigma)$ by
\begin{align*}
 \T(\sigma)=&\Lambda_{-\zeta}T^*X\otimes\bigotimes_{k=1}^\infty\left(
   \Lambda_{-\zeta q^k}T^*X\otimes \Lambda_{-\zeta^{-1}q^k}TX
   \otimes S_{q^k}T^*X\otimes S_{q^k}TX \right), \\
 \T_{\hatlam}(\sigma)=&
   \Lambda_{-\zeta}T^*\hatX_{\hatlam}\otimes
   \bigotimes_{k=1}^\infty\left(
   \Lambda_{-\zeta q^k}T^*\hatX_{\hatlam}\otimes 
   \Lambda_{-\zeta^{-1}q^k}T\hatX_{\hatlam}
   \otimes S_{q^k}T^*\hatX_{\hatlam}\otimes 
    S_{q^k}T\hatX_{\hatlam} \right), \\
 \W_{\hatlam,i}(\sigma)=&
   \Lambda_{-\zeta q^{f_{\hatlam,i}}}W^*_{\hatlam,i}\otimes
   \bigotimes_{k=1}^\infty\left(
   \Lambda_{-\zeta q^{f_{\hatlam,i}+k}}W^*_{\hatlam,i}\otimes 
   \Lambda_{-\zeta^{-1}q^{-f_{\hatlam,i}+k}}W_{\hatlam,i}\right) \\
   & \otimes S_{q^{f_{\hatlam,i}}}W^*_{\hatlam,i}\otimes
   \bigotimes_{k=1}^\infty\left(
   S_{q^{f_{\hatlam,i}+k}}W^*_{\hatlam,i}\otimes 
    S_{q^{-f_{\hatlam,i}+k}}W_{\hatlam,i} \right), \\
   \W_{\hatlam}(\sigma)=
   &\bigotimes_i\W_{\hatlam,i}(\sigma).
\end{align*}
Here $\Lambda_tW=\bigoplus_i \Lambda^iW$ and 
$S_tW=\bigoplus_iS^iW$ denote the total exterior power and total symmetric
power of a vector bundle $W$. We can write $\T$ and 
$\T_{\hatlam}\otimes\W_{\hatlam}$ 
as formal power series in $q$ and $q^{\frac{1}{r}}$
\begin{equation}\label{eq:Korb}
 \begin{split}
 \T(\sigma)&=\sum_{k=0}^\infty R_k(\sigma)q^k \\
 \T_{\hatlam}(\sigma)\otimes\W_{\hatlam}(\sigma)&
=\sum_{k=0}^\infty \hat{R}_{\hatlam,k}(\sigma)q^{\frac{k}{r}} 
 \end{split}
\end{equation}
with coefficients $R_k(\sigma),\hat{R}_{\hatlam,k}(\sigma)
\in K_{orb}(X)\otimes \Z[\zeta,\zeta^{-1}]$ 
where $K_{orb}(X)$ denotes the Grothendieck group 
of orbifold vector bundles and $r$ is the least common multiple 
of the orders $|H_\lam|$ of the isotropy groups. 

Let $W$ be a complex orbifold vector bundle over $X$ and 
$D\otimes W$ a spin-c Dirac operator twisted by $W$ 
in the sense of \cite{D}. 
It is an elliptic differential operator 
\begin{equation}\label{eq:Dirac} 
D\otimes W: \Gamma(X, E^+\otimes W)\to \Gamma(X,E^-\otimes W), 
\end{equation}
where 
\[ E^+=\bigoplus_{i:even}\Lambda^iTX \ \text{and}\ 
   E^-=\bigoplus_{i:odd}\Lambda^iTX. \]
$D\otimes W$ is constructed from hermitian 
metrics and connections on various vector bundles associated to
$TX$ and $W$. Its principal symbol is such that, when $X$ is a 
complex orbifold and $W$ is 
a holomorphic vector bundle, then $D\otimes W$ has the same principal symbol 
(up to multiplicative constant) 
as $\bar{\partial}+\bar{\partial}^*$ twisted by $W$ acting 
on the sections of $E^+\otimes W$. 

We now define
\begin{equation*}
 \begin{split}
 \var(X)=&\zeta^{-\frac{n}{2}}\ind(D\otimes \T(\sigma)) \\
 \hatvar(X)=&\zeta^{-\frac{n}{2}}\sum_{{\hatlam}\in \hatLam}
 \zeta^{f_{\hatlam}}\ind(D_{\hatX_{\hatlam}}\otimes\T_{\hatlam}(\sigma)
  \otimes \W_{\hatlam}(\sigma)),
\end{split} 
\end{equation*}
where $\dim X=2n$ and $D_{\hatX_{\hatlam}}$ is the spin-c Dirac operator 
for $\hatX_{\hatlam}$.
More precisely 
\begin{equation}\label{eq:var}
 \begin{split}
 \var(X)=&\zeta^{-\frac{n}{2}}\sum_{k=0}^\infty
 \ind(D\otimes R_k(\sigma))q^k, \\
 \hatvar(X)=&\zeta^{-\frac{n}{2}}\sum_{{\hatlam}\in \hatLam}
 \zeta^{f_{\hatlam}}\sum_{k=0}^\infty
 \ind(D_{\hatX_{\hatlam}}\otimes\hat{R}_{\hatlam,k}
(\sigma))q^{\frac{k}{r}}, 
\end{split} 
\end{equation}
where $R_k$ and $\hat{R}_k$ are given in \eqref{eq:Korb}. 
$\var(X)$ and $\hatvar(X)$ are called \emph{elliptic genus} and 
\emph{orbifold elliptic genus} of the orbifold $X$ respectively. 

Let $N>1$ be an integer. If $f=\frac{s}{r}$ is a rational number with 
$r$ relatively prime to $N$, let $d$ be an integer such that 
$dr\equiv 1 \bmod N$. 
We define $\brf\in \Z$ by
\[ \brf=ds .\]
$\brf$ is determined modulo $N$. It satisfies 
$\brf_1+\brf_2\equiv\breve{f_1+f_2}\ (\bmod\ N)$. 

Assume that $N$ is relatively prime 
to $|H_\lam|$ for all $\lam\in \Lambda$. Then $\brf_{\hatlam}$ is 
defined since $f_{\hatlam}$ can be written in the form 
$f_{\hatlam}=\frac{s}{|H_{\pi(\hatlam)}|}$. 
Under the above assumption we put $\sigma=\frac{k}{N}$ with 
$0<k<N$ and define the \emph{modified orbifold 
elliptic genus} $\brvar(X)$ of \emph{level $N$} by 
\begin{equation}\label{eq:brvar}
 \brvar(X)=\zeta^{-\frac{n}{2}}\sum_{{\hatlam}\in \hatLam}
  \zeta^{\brf_{\hatlam}}\ind(D_{\hatX_{\hatlam}}\otimes\T_{\hatlam} 
   \otimes\W_{\hatlam}).
\end{equation}

The genus $\zeta^{\frac{n}{2}}\var(X)$ belongs to 
$(\Z[\zeta,\zeta^{-1}])[[q]]$ and the genus 
$\zeta^{\frac{n}{2}}\hatvar(X)$ to 
$(\Z[\zeta^{\frac{1}{r}},\zeta^{-\frac{1}{r}}])[[q^{\frac{1}{r}}]]$. 
Similarly 
the genus $\zeta^{\frac{n}{2}}\brvar(X)$ belongs to 
$(\Z[\zeta,\zeta^{-1}]/(\zeta^N))[[q^{\frac{1}{r}}]]$. 
When it is necessary to make explicit the parameters $\tau$ and $\sigma$ 
we write $\var(X;\tau,\sigma), \hatvar(X;\tau,\sigma)$ and 
$\brvar(X;\tau,\sigma)$ for $\var(X),\hatvar(X)$ and $\brvar(X)$ respectively. 

When a compact connected Lie group $G$ acts on $X$ preserving 
almost complex structure, $G$ acts naturally on vector bundles 
$\T,\T_{\hatlam}$ and $\W_{\hatlam}$. This is clear for 
$\T$. As to the other two it should be noticed that 
the action of a connected group preserves each stratum and 
each conjugacy class of the isotropy group of the stratum. 
It follows that the action of $G$ lifts to the action on each sector, 
and hence on $\T_{\hatlam}$. Moreover some finite covering group of 
$G$ acts on each $W_{\hatlam}$ compatibly with the decomposition 
\eqref{eq:eigendecomp2} as we shall see in Section 4, Lemma \ref{lemm:liftable}. 
Then it also acts on $\W_{\hatlam}$. We shall assume here that 
$G$ itself acts on each $\W_{\hatlam}$. By using a $G$ invariant 
spin-c Dirac operator the formulae \eqref{eq:var} and 
\eqref{eq:brvar} define equivariant genera which will be denoted by the 
same symbols. In this case 
$\zeta^{\frac{n}{2}}\var(X)$ belongs to 
$(R(G)\otimes \Z[\zeta,\zeta^{-1}])[[q]]$, 
$\zeta^{\frac{n}{2}}\hatvar(X)$ to
$(R(G)\otimes\Z[\zeta^{\frac{1}{r}},\zeta^{-\frac{1}{r}}])[[q^{\frac{1}{r}}]]$
and $\zeta^{\frac{n}{2}}\brvar(X)$ to 
$(R(G)\otimes\Z[\zeta,\zeta^{-1}]/(\zeta^N))[[q^{\frac{1}{r}}]]$, 
where $R(G)$ is the character ring of the group $G$. 
The value at $g\in G$ of $\var(X)$ is denoted by 
$\var_g(X)$, and similarly by $\hatvar_g(X), \brvar_g(X)$ for
$\hatvar(X), \brvar(X)$.

Let $N>1$ be an integer. When $\sigma=\frac{k}{N},\ 0<k<N,$ the 
genera $\var(X)$ and $\hatvar(X)$ are also called of \emph{level} $N$. 
We can now state the main theorems of the present paper. 

\begin{thm}\label{thm:brvar}
Let $X$ be an almost complex closed orbifold of 
dimension $2n$ with a non-trivial 
action of a compact connected Lie group $G$. Let $N>1$ be an integer 
such that $N$ is relatively prime to the orders of all isotropy 
groups $H_\lambda$. Assume that there 
is an orbifold line bundle $L$ with a lifted action of $G$ over $X$ 
such that $\Lambda^nTX=L^N$. Then the equivariant 
modified elliptic genus 
$\brvar_g(X)$ of level $N$ is rigid, that is, $\brvar_g(X)$ is constant as 
a function on $G$ for $\sigma=\frac{k}{N},\ 0<k<N.$ 
\end{thm}

\begin{rem}
The Picard group of orbifold line bundles over $X$ is isomorphic 
to the second cohomology group $H^2(X,\Z_X)$ where $\Z_X$ is a certain 
sheaf over $X$, cf. \cite{SW}. This correspondence can be considered as 
assigning the first Chern class to orbifold line bundles. In this sense 
there is an orbifold line bundle $L$ satisfying the condition 
$\Lambda^nTX=L^N$ 
if and only if the first Chern class $c_1(\Lambda^nTX)$ is divisible by 
$N$ in $H^2(X,\Z_X)$. The class $c_1(\Lambda^nTX)$ can 
be called the first Chern class of the 
almost complex orbifold $X$ and be written $c_1(X)$ as the manifold case. 
\end{rem}

Let $L\to X$ be an orbifold line bundle and let 
$\tilde{\pi}:(\tilde{V},\tilde{U},H)\to(V,U,H)$ be an orbifold chart 
of $L$. $L$ will be called a \emph{genuine} line 
bundle if $h$ acts trivially on the fiber 
$\tilde{\pi}^{-1}(v)$ over $v$ for any $h$ such that $hv=v$. In this
case $\pi: L\to X$ becomes a line bundle over the space $X$ 
in the usual sense. The orbifold elliptic genus is not rigid 
in general. However we have the following theorem.

\begin{thm}\label{thm:hatvar}
Let $X$ be an almost complex closed orbifold of 
dimension $2n$ with a non-trivial 
action of a compact connected Lie group $G$. Let $N>1$ be an integer. 
Assume that there 
is a genuine line bundle $L$ with a lifted action of $G$ over $X$ such that
$\Lambda^nTX=L^N$. Then the equivariant elliptic genus 
$\hatvar_g(X)$ of level $N$ is rigid, that is, $\hatvar_g(X)$ is constant as 
a function on $G$ for $\sigma=\frac{k}{N},\ 0<k<N.$ 
\end{thm}
\begin{note}
Under the assumption of Theorem \ref{thm:hatvar} each $f_{\hatlam}$ 
is an integer. See Note after Lemma \ref{lemm:hatvarA}. 
\end{note}

The next Theorem concerns the case where $\Lambda^nTX$ is a torsion 
element in the Picard group of orbifold line bundles.
\begin{thm}\label{thm:hatvartorsion}
Let $X$ be an almost complex closed orbifold of 
dimension $2n$ with a non-trivial 
action of a compact connected Lie group $G$. 
Assume that $\Lambda^nTX$ is trivial as an orbifold line bundle. 
Then the equivariant elliptic genus 
$\hatvar_g(X)$ is rigid, that is, $\hatvar_g(X)$ is constant as 
a function on $G$.
\end{thm}

Theorem \ref{thm:hatvartorsion} is essentially due to Dong, Liu 
and Ma \cite{DLM} in a more general setting. 
It seems the hypothesis concerning the vanishing of first 
Chern class in \cite{DLM} was erroneously stated. 

The proof of Theorem \ref{thm:brvar}, Theorem \ref{thm:hatvar} and 
Theorem \ref{thm:hatvartorsion} 
will be given in Section 5 and Section 6.

\section{Vergne's fixed point formula}

In this section we review the fixed point formula due to Vergne 
\cite{V}. 

Let $X$ be an almost complex closed orbifold 
of dimension $2n$ 
and $\pi:W\to X$ a complex orbifold vector bundle. Introduce 
a hermitian metric and a hermitian connection on $W$, and 
let $\omega$ be the curvature form of the connection.
It is a differential form with values in $End(W)$ where $End(W)$ is 
the vector bundle of skew hermitian endomorphisms of $W$. 
Locally it is an $H_x$-invariant 
differential form $\omega_V$ on $V$ with values in 
$End(\tilde{V})$, where $(V,U,H)$ is a chart 
for $X$ and $(\tilde{V},\tilde{U},H)$ is a chart for
the bundle $W\to X$. These $\omega_V$ behave compatibly with 
injections of charts. With respect to an 
orthonormal basis $\omega_V$ is expressed as a skew-hermitian matrix valued 
form.

The form $\Gamma(W)=-\frac{1}{2\pi\img}\omega$ is called 
the \emph{Chern matrix}, cf. \cite{D}. The form 
\[ c(W)=\sum_{i=1}^{d}c_i(W)=\det(1+\Gamma(W))\quad (d=\rank W) \] 
is the total Chern form of $W$. 
It is convenient to write formally 
\begin{equation}\label{eq:split}
 c(W)=\prod_{i=1}^{d}(1+x_i). 
\end{equation}
The $x_i$ are called the Chern roots and \eqref{eq:split} is 
called the formal splitting of the total Chern form. 
The Chern character of $W$ is 
defined by 
\[ ch(W)=\tr e^{\Gamma(W)}=\sum_{i=1}^{d}e^{x_i}. \]

The Todd form of the almost complex orbifold $X$ is given by 
\[ Td(X)=\det(\frac{\Gamma(TX)}{1-e^{-\Gamma(TX)}})=
 \prod_{i=1}^{n}\frac{x_i}{1-e^{-x_i}}, \]
with $c(TX)=\prod_i(1+x_i)$. 

Now suppose that a compact connected Lie group $G$ acts on $X$ preserving 
almost complex structure. Let $W$ be a complex orbifold 
vector bundle over $X$ with a compatible action of $G$. 
Let $F=X^G$ be the fixed point set of the action. 
$F$ is an almost complex orbifold and the 
normal bundle of each component of $F$ is a complex orbifold 
vector bundle. 
We can say more about this. 

\begin{lemm}\label{lemm:liftable}
Let $(V_x,U_x,H_x)$ be a reduced chart centered at $x\in F$
such that $U_x$ is invariant under the action of $G$. 
Then the action lifts to an action of some finite covering group
$\tilde{G}\to G$ on $V_x$. That action commutes with the action of $H_x$. 
\end{lemm}
\begin{proof}
From the fact that $G$ acts 
on $X$ it follows that there are an automorphism $\rho$ of 
the group $H_x$, a small neighborhood $O$ of the identity element 
$e$ of $G$, and a local transformation group 
$\psi:O\times V_x\to V_x $ such that 
\[ \psi_g\circ h=\rho(h)\circ \psi_g, \]
where $\psi_g:V_x\to V_x$ is defined by $\psi_g(v)=\psi(g,v)$. 
We may assume that $\psi_e$ is the identity map by conjugating 
by an element of $H_x$ if necessary. It follows that $\rho$ is the identity 
automorphism, and hence $\psi_g$ commutes with the action of $H_x$. 
The local transformation group $\psi$ extends to an action of 
the universal covering group $G^{univ}$ of $G$. 
Let $H_0$ be the normal subgroup 
of $H_x$ consisting of elements which act trivially on $V_x$. The 
quotient $H_x/H_0$ acts effectively on $V_x$ and if 
$v$ is a generic point in $V_x$, then $H_x/H_0$ acts simply transitively 
on the orbit of $v$. We identify $\pi_1(G)$ with the kernel of 
the projection $G^{univ}\to G$. Then $\pi_1(G)$ acts through $\psi$ 
on each orbit of $H_x$ and the action $\psi$ commutes 
with that of $H_x$. Hence we get a homomorphism 
$\alpha:\pi_1(G)\to H_x/H_0$ defined by 
\[ \alpha(g)v=\psi_g(v). \] 
Let $\Gamma$ be the kernel of $\alpha$. Then $\tilde{G}=G^{univ}/\Gamma$ is a 
finite covering of $G=G^{univ}/\Gamma$ and the action $\psi_g$ induces the 
action of $\tilde{G}$ on $V_x$. 
\end{proof}
Let $F_\nu$ be a component of $F$ and $x\in F_\nu$. If a finite covering 
group
$\tilde{G}_\nu$ of $G$ acts on $V_x$, then $\tilde{G}_\nu$ also acts on 
$V_y$ for any $y\in U_x$. It follows that $\tilde{G}_\nu$ acts 
on the normal bundle $N(F_\nu,X)$ of $F_\nu$ in $X$ by (fiberwise) 
automorphisms. Since there are only a finite number of components in $F$ 
and there is a common finite covering group $\tilde{G}$ of all the 
$\tilde{G}_\nu$, it follows that the group $\tilde{G}$ acts on each 
$N(F_\nu,X)$. Similar consideration yields that we can take a 
finite covering group $\tilde{G}$ which acts on the normal bundle 
$W_{\hatlam}$ of each sector $X_{\hatlam}$. 
Hereafter we assume that \emph{$G$ itself acts on $V_x$} for any $x\in F$ 
replacing $G$ by $\tilde{G}$ if necessary. Then $G$ acts also 
on the normal bundle $N(F,X)$ of $F$ in $X$ by (fiberwise) automorphisms.   

For the fixed point formula it is enough to assume that the group $G$ is 
a torus $T$. We then take a topological generator $g\in T$. 
Let $\{\hatF_{\hatlam}\}_{\hatlam\in \hatLam_F}$ be the 
totality of sectors of $F=X^T$. Charts of $\hatF_{\hatlam}$ 
are of the form
\[ (V_x^{g,h}, V_x^{g,h}/C(h),C(h)), \]
where $V_x^{g,h}=(V_x^g)^h=(V_x^h)^g$. 
The normal bundle of the immersion $\pi :\hatF_{\hatlam}\to F$ is 
denoted by $N_{\hatlam}$. We also set 
$\hatN_{\hatlam}=\pi^*N(F,X)\vert \hatF_{\hatlam}$ 
and $\hatW_{\hatlam}=\pi^*W\vert \hatF_{\hatlam}$. 
If $\gamma_x=[h]\in \calC(H_x)$, then 
$h$ acts on each fiber of 
$N_{\hatlam},\hatN_{\hatlam}$ and $\hatW_{\hatlam}$ 
in the sense as explained in Section 3. 
Also $T$ acts on $\hatN_{\hatlam}$ and $\hatW_{\hatlam}$. 

We define three differential forms on $\hatF_{\hatlam}$ by 
\begin{equation*}\label{eq:chform}
 \begin{split}
 D_h(N_{\hatlam})&=\det(1-h^{-1}e^{-\Gamma(N_{\hatlam})}) \\
 D_{g,h}(\hatN_{\hatlam})&=
      \det(1-g^{-1}h^{-1}e^{-\Gamma(\hatN_{\hatlam})}) \\
 ch_{g,h}(\hatW_{\hatlam})&=\tr(ghe^{\Gamma(\hatW_{\hatlam})}).
 \end{split}
\end{equation*}

Vergne's fixed point formula is stated in the following form. 
Let $D$ be a spin-c Dirac operator on an almost complex 
orbifold $X$. Then we have 
\begin{equation}\label{eq:fixpoint}
 \ind_t(D\otimes W)=\sum_{\hatlam\in \hatLam_F}
 \frac{1}{m(\hatF_{\hatlam})}\int_{\hatF_{\hatlam}}
 \frac{Td(\hatF_{\hatlam})\ch_{g,h}(\hatW_{\hatlam})}
      {D_h(N_{\hatlam})D_{g,h}(\hatN_{\hatlam})}.
\end{equation}
Note that $m(\hatF_{\hatlam})=|C(h)|$ where $[h]=\gamma_x$ is 
an arbitrary point in $\hatF_{\hatlam}^{gen}$. Let $h'$ be another 
representative of $[h]$. Then the identification 
$V_x^h/C(h)=V_x^{h'}/C(h')$ induces the equality 
\[ \frac{\ch_{g,h}(\hatW_{\hatlam})}
     {D_h(N_{\hatlam})D_{g,h}(\hatN_{\hatlam})}
 =  \frac{\ch_{g,h'}(\hatW_{\hatlam})}
      {D_{h'}(N_{\hatlam})D_{g,h'}(\hatN_{\hatlam})}. \]
Since $|\gamma_x||C(h)|=|H_x|=|H_{\pi(\hatlam)}|$ we obtain 
from \eqref{eq:fixpoint} 
\begin{equation}\label{eq:fixpoint2}
 \ind_t(D\otimes W)=\sum_{\hatlam\in \hatLam_F}
 \frac{1}{|H_{\pi(\hatlam)}|}\int_{\hatF_{\hatlam}}\sum_{h\in \gamma_x}
 \frac{Td(\hatF_{\hatlam})\ch_{g,h}(\hatW_{\hatlam})}
      {D_h(N_{\hatlam})D_{g,h}(\hatN_{\hatlam})},
\end{equation}
where $\gamma_x\in \hatF_{\hatlam}^{gen}$. Since 
$\hatF_{\hatlam}\setminus\hatF_{\hatlam}^{gen}$ has at least codimension 
two, the above integral is well-defined. 

The following observation is useful later. Suppose that $T=S^1$. 
We put $t=e^{2\pi\img z}$ with $z\in \R$.  
$\hatN_{\hatlam}$ decomposes into the direct sum
\[\hatN_{\hatlam} =\bigoplus_{i,j}E_{\chi_i^{S^1},\chi_j^{\l h\r}}\quad  
\text{(finite sum)} \]
of eigen-bundle where the sum is extended over the pairs of characters
$(\chi_i^{S^1},\chi_j^{\l h\r})$ of $S^1$ and $\l h\r
(=\text{cyclic subgroup of $H_x$ generated by $h$)}$, and $(t,h)$ acts on 
$E_{\chi_i^{S^1},\chi_j^{\l h\r}}$ by multiplication by 
$\chi_i^{S^1}(t)\chi_j^{\l h\r}(h)$. Write 
\[ \chi_i^{S^1}(t)=e^{2\pi\img m_iz} \ \text{with $t=e^{2\pi\img z}$ and}\ 
\chi_j^{\l h\r}(h)=e^{2\pi\img m_j(h)}, \]
and $c(E_{\chi_i^{S^1},\chi_j^{\l h\r}})=\prod(1+x_{i,j,k})$ formally 
as before. We put $x_{i,j,k}=2\pi\img y_{i,j,k}$. Then 
\[ D_{t,h}(\hatN_{\hatlam})=\prod_{i,j}
   D_{t,h}(E_{\chi_i^{S^1},\chi_j^{\l h\r}}) \]
with
\begin{equation}\label{eq:Dth}
 D_{t,h}(E_{\chi_i^{S^1},\chi_j^{\l h\r}})= 
 \det(1-t^{-1}h^{-1}e^{-\Gamma(E_{\chi_i^{S^1},\chi_j^{\l h\r}})})=
 \prod_k (1-e^{2\pi\img(-m_iz-m_j(h)-y_{i,j,k})}). 
\end{equation}
Similar formulae for $D_h(N_{\hatlam})$ and 
$ch_{t,h}(\hatW_{\hatlam})$ are available.

\section{Proof of main theorems; I fixed point formula}

In order to prove Theorem \ref{thm:brvar}, Theorem\ref{thm:hatvar} 
and Theorem \ref{thm:hatvartorsion} it is enough to assume that $G$ is a 
circle group $S^1$ since the character of a compact connected 
Lie group $G$ is determined by those of all circle subgroups of $G$. 
Let $X$ be a connected stably almost complex closed orbifold 
with a non-trivial $S^1$ action. 
Let $\hatX_{\hatlam}$ be a sector of $X$ where 
$\hatlam\in \hatLam$ as in Section 2. As remarked in Section 3, 
$S^1$ may be assumed to act on $\hatX_{\hatlam}$. 
We first investigate sectors of the orbifold $\hatX_{\hatlam}^{S^1}$. 
For that purpose we need some notations. 
Let $H$ be a finite group. We define 
\[ CM(H)=\{(h_1,h_2)\in H\times H \mid h_1h_2=h_2h_1\}. \]
The group acts on $CM(H)$ by conjugations on each factor. Let
$\hatcalC(H)$ be the set of conjugacy classes by that action. 
If $\delta=[h_1,h_2]\in \hatcalC(H)$ is the conjugacy class of 
$(h_1,h_2)$, we define $\delta^{(i)}\in \calC(H),\ i=1,2$, by
$\delta^{(1)}=[h_1],\ \delta^{(2)}=[h_2]$. A group homomorphism 
$H\to H'$ induces an obvious map $\hatcalC(H)\to \hatcalC(H')$.

Let $\hatpi:\HatF\to \hatF=\hatX^{S^1}$ be the total sector. 
\begin{lemm}\label{lemm:Hat}
The inverse image $(\pi\circ\hatpi)^{-1}(x)$ for $x\in F=X^{S^1}$ 
can be identified with $\hatcalC(H_x)$ in a canonical way. 
With this identification made, an element 
$\delta\in \hatcalC(H_x)$ belongs to a sector of $\hatX_{\hatlam}$ if 
and only if $\delta^{(1)}\in \hatX_{\hatlam}$.   
\end{lemm}
\begin{proof}
Take an element $\hatgam_x\in (\pi\circ\hatpi)^{-1}(x)$ and 
put $\gamma_x=\hatpi(\hatgam_x)$. If $h_1\in H_x$ is a representative 
of $\gamma_x$ then $\hatgam_x$ is identified with a conjugacy 
class $[h_2]$ of the group $C(h_1)$. 
We assign the equivalence class $[h_1,h_2]\in\hatcalC(H_x)$ of 
$(h_1,h_2)$ to $\hatgam_x$. This assignment is well-defined and 
bijective as is easily seen. 

If $\delta$ corresponds to $\hatgam_x$ by this identification, then 
$\delta^{(1)}=\gamma_x$. Hence $\hatpi(\delta)\in \hatX_{\hatlam}$ 
if and only if $\delta^{(1)}\in \hatX_{\hatlam}$.
\end{proof}

Let $\mathcal{H}$ be the upper half plane. 
We consider the function $\Phi(z,\tau)$ defined on 
$\C\times \mathcal{Hi}$ 
given by the following formula. 
\[ \Phi(z,\tau)=(t^{\frac12}-t^{-\frac12})\prod_{k=1}^{\infty}
 \frac{(1-tq^k)(1-t^{-1}q^k)}
      {(1-q^k)^2},  \]
where $t=e^{2\pi\img z}$ and $q=e^{2\pi\img\tau}$. Note that $|q|<1$. 

The group $SL_2(\Z)$ 
acts on $\C\times\mathcal{H}$ by 
\[ A(z,\tau)=(\frac{z}{c\tau+d},A\tau)
 =(\frac{z}{c\tau+d}, \frac{a\tau+b}{c\tau+d}),\quad 
A=\begin{pmatrix} a & b \\ c & d \end{pmatrix} .\]

$\Phi$ is a Jacobi form and satisfies the following transformation
formulae, cf. \cite{HBJ}. 
\begin{equation}\label{eq:transformula1}
 \begin{split}
 \Phi(A(z,\tau))&=(c\tau +d)^{-1}
  e^{\frac{\pi\img cz^2}{c\tau +d}}\Phi(z,\tau), \\
 \Phi(z+m\tau +n,\tau)&=(-1)^{m+n}e^{-\pi\img(2mz+m^2\tau)}\Phi(z,\tau)
 \end{split}
\end{equation}
where $m,n\in \Z$.

For $\sigma\in \C$ we set 
\[ \phi(z,\tau,\sigma)=\frac{\Phi(z+\sigma,\tau)}{\Phi(z,\tau)}
=\zeta^{-\frac12}\frac{1-\zeta t}{1-t}
\prod_{k=1}^{\infty}\frac{(1-\zeta tq^k)(1-\zeta^{-1}t^{-1}q^k)}
 {(1-tq^k)(1-t^{-1}q^k)}, \]
where $\zeta =e^{2\pi\img\sigma}$. 
From \eqref{eq:transformula1} 
the following transformation formulae 
for $\phi$ follow:
\begin{equation}\label{eq:transformula2}
 \begin{split}
 \phi(A(z,\tau),\sigma)&=
  e^{\pi\img c(2z\sigma +(c\tau +d)\sigma^2)}\phi(z,\tau,(c\tau +d)\sigma), \\
  \phi(z+m\tau +n,\tau,\sigma)&=e^{-2\pi\img m\sigma}\phi(z,\tau,\sigma)=
  \zeta^{-m}\phi(z,\tau,\sigma).
 \end{split}
\end{equation}

For the later use we extend the domain of the function $\phi$. 
Let $w$ be a polynomial in an indeterminate $y$. We put 
\[ \phi(w,\tau,\sigma)=\zeta^{-\frac12}\frac{1-\zeta e^{2\pi\img w}}
 {1-e^{2\pi\img w}}
\prod_{k=1}^{\infty}\frac{(1-\zeta e^{2\pi\img w}q^k)
  (1-\zeta^{-1}e^{-2\pi\img w}q^k)}
 {(1-e^{2\pi\img w}q^k)(1-e^{-2\pi\img w}q^k)}. \]
and consider it formally as an element of $\C[[y]]$. 
In the above extended meaning the function $\phi$ still 
satisfies the same transformation laws:
\begin{equation}\label{eq:transformula3}
\begin{split}
 \phi(A(w,\tau),\sigma)&=
  e^{\pi\img c(2w\sigma +(c\tau +d)\sigma^2)}\phi(w,\tau,(c\tau +d)\sigma), \\
  \phi(w+m\tau +n,\tau,\sigma)&=e^{-2\pi\img m\sigma}\phi(w,\tau,\sigma)=
  \zeta^{-m}\phi(w,\tau,\sigma).
 \end{split}
\end{equation}
In fact, if we substitute an arbitrary complex number for $y$ 
then the above equalities hold by \eqref{eq:transformula2}. Hence 
the equalities \eqref{eq:transformula3} hold in $\C[[y]]$. 
If $A^*=\sum_{l=0}^\infty A^{2l}$ is a commutative graded algebra over $\C$ 
with even grading and with $A^0=\C$, then $\phi(w,\tau,\sigma)$ 
can be defined for $w\in A^*$ and it satisfies \eqref{eq:transformula3} 
since one may think of $w$ as a specilization of an element in $\C[y]$. 

Hereafter we assume that $X$ is an almost complex closed orbifold of 
dimension $2n$ with a non-trivial action of $S^1$. Let $\HatF_{\Hatlam}$ be 
a sector of $\hatF=\hatX^{S^1}$, that is a component of $\HatF$. 
Its orbifold charts are of the form 
\[ (V_x^{t,h_1,h_2}, V_x^{t,h_1,h_2}/(C(h_1)\cap C(h_2)), C(h_1)\cap C(h_2)), 
(h_1,h_2)\in CM(H_x). \] 
Note that $C(h_1)\cap C(h_2)$ coincides with $C_{C(h_1)}(h_2)$, the 
centralizer of $h_2$ in $C(h_1)$. 
The image $\hatpi(\HatF_{\Hatlam})$ is contained in a unique 
sector $\hatX_{\hatlam}$ of the orbifold $X$. 
We formally write
\begin{equation}\label{eq:split2}
 \begin{split}
 c(\HatF_{\Hatlam})&=\prod_{i=1}^{r_1}(1+x_i), \\
 c(N_{\Hatlam})&=\prod_{i=r_1+1}^{r_1+r_2}(1+x_i), \\ 
c(\HatN_{\Hatlam})&=\prod_{i=r_1+r_2+1}^{r_1+r_2+r_3}(1+x_i), \\ 
c(\hatW_{\Hatlam})&=\prod_{i=r_1+r_2+r_3+1}^{r_1+r_2+r_3+r_4}(1+x_i), 
 \end{split}
\end{equation}
where $N_{\Hatlam}$ is the normal bundle of the immersion 
$\HatF_{\Hatlam}\to \hatX_{\hatlam}^{S^1},\  
\HatN_{\Hatlam}=\hatpi^{-1}N(\hatX_{\hatlam}^{S^1},
 \hatX_{\hatlam})|\HatF_{\Hatlam}$ 
and $\hatW_{\Hatlam}=\hatpi^{-1}W_\lam|\HatF_{\Hatlam}$. Recall that 
$W_\lam$ is the normal bundle of the immersion $\hatX_{\hatlam}\to X$, 
cf. Section 3. 
Note that $r_1+r_2+r_3+r_4=n$. We denote the Euler class of the orbifold
$\HatF_{\Hatlam}$ by $e(\HatF_{\Hatlam})$. It is equal to the top 
Chern class $c_{r_1}(\HatF_{\Hatlam})$ and is written 
\[ e(\HatF_{\Hatlam})=\prod_{i=1}^{r_1}x_i .\]
Finally we put 
\[ x_i=2\pi\img y_i .\]

Take a point $\hatgam_x\in \HatF_{\Hatlam}$ and identify it 
with $\delta=[h_1,h_2]\in \hatcalC(H_x)$. Then $h_1$ acts on $W_\lam$ and 
hence on $\hatW_{\Hatlam}$ as was explained in Section 3. 
Also as was explained in Section 4, $h_2$ acts on 
$N_{\Hatlam}, \HatN_{\Hatlam}$ and $\hatW_{\Hatlam}$. Furthermore 
$t=e^{2\pi\img z}\in S^1$ acts on $\HatN_{\Hatlam}$ and $\hatW_{\Hatlam}$. 
The weights of these actions can be taken 
compatibly with the formal splitting \eqref{eq:split2}. We write them 
in the following form
\begin{equation}\label{eq:weights}
 \begin{cases}
 m_i^{S^1}  &\text{for $t=e^{2\pi\img z}$}, \\
 m_i^{(h_1,h_2)}(h_1)\ &\text{for $h_1$}, \\
 m_i^{(h_1,h_2)}(h_2)\ &\text{for $h_2$}
\end{cases}
\end{equation}
for $i=1,\ldots ,r_1+r_2+r_3+r_4$. 
Note that $m_i^{S^1}\in \Z$ and $m_i^{(h_1,h_2)}(h_j),\ j=1,2,$ is 
a rational number determined modulo $\Z$. In the sequel we shall fix 
one representative $m_i^{(h_1,h_2)}(h_j),\ j=1,2,$ for each 
pair $(h_1,h_2)\in CM(H_x)$. 
We make the convention that 
\begin{equation}\label{eq:weights2} 
 \begin{split}
 m_i^{S^1}&=0 \quad \text{for $1\leq i\leq r_1+r_2$}, \\
 m_i^{(h_1,h_2)}(h_1)&=0 \quad \text{for $1\leq i\leq r_1+r_2+r_3$}, \\
 m_i^{(h_1,h_2)}(h_2)&=0 \quad \text{for $1\leq i\leq r_1$}.
 \end{split} 
\end{equation}

Hereafter we shall write $\hatvar(X;z,\tau,\sigma)$ and 
$\brvar(X;z,\tau,\sigma)$ instead of $\hatvar_t(X;\tau,\sigma)$ and 
$\brvar_t(X;\tau,\sigma)$. Let $\{\HatF_{\Hatlam}\}_{\Hatlam\in \HatLam}$ be 
the totality of sectors of $\hatF$. 

\begin{prop}\label{prop:fxtptformula}
Let $X$ be an almost complex orbifold with a non-triavial action of $S^1$. 
Then the equivariant elliptic orbifold elliptic genus 
$\hatvar(X;z,\tau,\sigma)$ is given by 
\begin{equation}\label{eq:fxtptformula}
\begin{split} 
 \hatvar(X;z,\tau,\sigma &)=\sum_{\Hatlam\in\HatLam}
 \frac{1}{|H_{\pi\circ\hatpi(\delta)}|} 
 \int_{\HatF_{\Hatlam}}\sum_{(h_1,h_2)\in \delta}
 e(\HatF_{\Hatlam})\cdot \\
 & \prod_{i=1}^ne^{2\pi\img m_i^{(h_1,h_2)}(h_1)\sigma}
 \phi(-y_i-m_i^{S^1}z+m_i^{(h_1,h_2)}(h_1)\tau-m_i^{(h_1,h_2)}(h_2),
  \tau,\sigma).
 \end{split} 
\end{equation}
Here $\delta$ is a point in $\HatF_{\Hatlam}$ as in 
Lemma \ref{lemm:Hat}. 
\end{prop}
\begin{note}
The above expressions give well-defined functions independent of 
the choice of representatives
$m_i^{(h_1,h_2)}(h_1),m_i^{(h_1,h_2)}(h_2)$ 
as is easily seen from \eqref{eq:transformula2}.
They are meromorphic functions in the variables $z,\tau,\sigma$. 
\end{note}

\begin{prop}\label{prop:brfxtptformula}
Let $N>1$ be an integer. We assume that $|H_x|$ is relatively prime 
to $N$ for all $x\in X$. Then the modified 
orbifold elliptic genus $\brvar(X;z,\tau,\sigma)$ of level 
$N$ is given by 
\begin{equation}\label{eq:brfxtptformula}
\begin{split} 
 \brvar(X;z,\tau,\sigma &)=\sum_{\Hatlam\in\HatLam}
 \frac{1}{|H_{\pi\circ\hatpi(\delta)}|} 
 \int_{\HatF_{\Hatlam}}\sum_{(h_1,h_2)\in \delta}
 e(\HatF_{\Hatlam})\cdot \\
 & \prod_{i=1}^ne^{2\pi\img \breve{m}_i^{(h_1,h_2)}(h_1)\sigma}
 \phi(-y_i-m_i^{S^1}z+m_i^{(h_1,h_2)}(h_1)\tau-m_i^{(h_1,h_2)}(h_2),
  \tau,\sigma).
 \end{split} 
\end{equation} 
\end{prop}
\begin{proof}
Let $\hatX_{\hatlam}$ be a sector of $X$. Recall that the contribution 
to $\hatvar(X)$ from $\hatX_{\hatlam}$ is 
\[ \zeta^{-\frac{1}{2}}\zeta^{f_{\hatlam}}
  \ind(D_{\hatX_{\hatlam}}\otimes\T_{\hatlam}\otimes\W_{\hatlam}). \]

If $\HatF_{\Hatlam}$ is a sector of $\hatX_{\hatlam}$ and 
$\delta\in \HatF_{\Hatlam}$, then $\delta^{(1)}$ lies in 
$\hatX_{\hatlam}$. For a moment we fix a representative 
$h_1$ of $\delta^{(1)}$. Then $\delta^{(2)}$ is a conjugacy 
class of $C(h_1)$. We apply the fixed point formula 
\eqref{eq:fixpoint2} to this and get 
\[
 \begin{split}
 \zeta^{-\frac{n}{2}}\zeta^{f_{\hatlam}}
 & \ind(D_{\hatX_{\hatlam}}\otimes\T_{\hatlam}\otimes\W_{\hatlam})= \\
 &\zeta^{-\frac{n}{2}}\sum_{\Hatlam\in \HatLam, \hatpi(\Hatlam)=\hatlam}
 \frac{1}{|C(h_1)|}
 \zeta^{f_{\hatlam}}\int_{\HatF_{\Hatlam}}\sum_{h_2\in \delta^{(2)}}
 \frac{Td(\HatF_{\Hatlam})ch_{t,h_2}(\T_{\Hatlam})
 ch_{t,h_1,h_2}(\W_{\Hatlam})}
 {D_{h_2}(N_{\Hatlam})D_{t,h_2}(\HatN_{\Hatlam})}
\end{split}
\]
where $\T_{\Hatlam}=\hatpi^*(\T_{\hatlam})|\HatF_{\Hatlam}$ and 
$\W_{\Hatlam}=\hatpi^*(\W_{\hatlam})|\HatF_{\Hatlam}$. Note that
$\hat{\pi}^*(T\hatX_{\hatlam})|\HatF_{\Hatlam}=
T\HatF_{\Hatlam}\oplus N_{\Hatlam}
 \oplus \HatN_{\Hatlam}$.
Using \eqref{eq:Dth}, \eqref{eq:split2} and \eqref{eq:weights} we have 
\begin{align*}
 Td(\HatF_{\Hatlam})&=\prod_{i=1}^{r_1}\frac{x_i}{1-e^{-2\pi\img y_i}}, \\
 D_{h_2}(N_{\Hatlam})&=\prod_{i=r_1+1}^{r_1+r_2}
  (1-e^{-2\pi\img(y_i+m_i^{(h_1,h_2)}(h_2))}), \\
 D_{t,h_2}(\HatN_{\Hatlam})&= \prod_{i=r_1+r_2+1}^{r_1+r_2+r_3}
 (1-e^{-2\pi\img(y_i+m_i^{S^1}z+m_i^{(h_1,h_2)}(h_2))}). 
\end{align*}
We also have
\[ \begin{split}
 ch_{t,h_2}&(\T_{\Hatlam})=\prod_{i=1}^{r_1+r_2+r_3}\left(
(1-\zeta e^{-2\pi\img(y_i+m_i^{S^1}z+m_i^{(h_1,h_2)}(h_2))})\cdot\right. \\
 &\left.\prod_{k=1}^\infty
 \frac{(1-\zeta q^ke^{-2\pi\img(y_i+m_i^{S^1}z+m_i^{(h_1,h_2)}(h_2))})
 (1-\zeta^{-1} q^k e^{2\pi\img(y_i+m_i^{S^1}z+m_i^{(h_1,h_2)}(h_2))})}
 {(1-q^ke^{-2\pi\img(y_i+m_i^{S^1}z+m_i^{(h_1,h_2)}(h_2))})
 (1-q^k e^{2\pi\img(y_i+m_i^{S^1}z+m_i^{(h_1,h_2)}(h_2))})}\right),
 \end{split} \]
and  
\[ \begin{split}
&ch_{t,h_1,h_2}(\W_{\Hatlam})=\prod_{i=r_1+r_2+r_3+1}^{r_1+r_2+r_3+r_4}
 \left(\frac{(1-\zeta q^{f_{\hatlam,i}}
 e^{-2\pi\img(y_i+m_i^{S^1}z+m_i^{(h_1,h_2)}(h_2))})}
 {(1-q^{f_{\hatlam,i}}
 e^{-2\pi\img(y_i+m_i^{S^1}z+m_i^{(h_1,h_2)}(h_2))})}\cdot\right. \\ 
 &\left.\prod_{k=1}^\infty
 \frac{(1-\zeta q^{f_{\hatlam,i}+k}e^{-2\pi\img(y_i+m_i^{S^1}z+
 m_i^{(h_1,h_2)}(h_2))})
 (1-\zeta^{-1} q^{-f_{\hatlam,i}+k} e^{2\pi\img(y_i+m_i^{S^1}z+
 m_i^{(h_1,h_2)}(h_2))})}
 {(1-q^{f_{\hatlam,i}+k}e^{-2\pi\img(y_i+m_i^{S^1}z+
 m_i^{(h_1,h_2)}(h_2))})
  (1-q^{-f_{\hatlam,i}+k} e^{2\pi\img(y_i+m_i^{S^1}z+
 m_i^{(h_1,h_2)}(h_2))})}\right).
\end{split} \]
Gathering these together we obtain
\begin{equation}\label{eq:hatvarf} 
\begin{split}
\zeta^{-\frac{n}{2}}\zeta^{f_{\hatlam}}
  \ind(D_{\hatX_{\hatlam}}&\otimes\T_{\hatlam}\otimes\W_{\hatlam})=  
 \sum_{\Hatlam\in \HatLam, \hatpi(\Hatlam)=\hatlam}
 \frac{1}{|C(h_1)|}\int_{\HatF_{\Hatlam}}\sum_{h_2\in \delta^{(2)}}
 e(\HatF_{\Hatlam})\cdot \\
 &\prod_{i=1}^ne^{2\pi\img f_{\hatlam,i}\sigma}
 \phi(-y_i-m_i^{S^1}z+f_{\hatlam,i}\tau-m_i^{(h_1,h_2)}(h_2),
  \tau,\sigma). 
 \end{split} 
\end{equation}
Note that $f_{\hatlam,i}\equiv m_i^{(h_1,h_2)}(h_1) \mod \Z$ for 
$r_1+r_2+r_3+1\leq i\leq r_1+r_2+r_3+r_4$ and $f_{\hatlam,i}=0$ 
for $i\leq r_1+r_2+r_3$. 
Then, by \eqref{eq:transformula3},
\[ \begin{split}
 e^{2\pi\img f_{\hatlam,i}\sigma}
 &\phi(-y_i-m_i^{S^1}z+f_{\hatlam,i}\tau-m_i^{(h_1,h_2)}(h_2),
  \tau,\sigma)= \\
 &e^{2\pi\img m_i^{(h_1,h_2)}(h_1)\sigma}
 \phi(-y_i-m_i^{S^1}z+m_i^{(h_1,h_2)}(h_1)\tau-m_i^{(h_1,h_2)}(h_2),
  \tau,\sigma).
\end{split} \] 

Putting this into \eqref{eq:hatvarf} we get 
\begin{equation*}\label{eq:hatvarm} 
\begin{split}
\zeta^{-\frac{1}{2}}\zeta^{f_{\hatlam}}&
  \ind(D_{\hatX_{\hatlam}}\otimes\T_{\hatlam}\otimes\W_{\hatlam})= 
 \sum_{\Hatlam\in \HatLam, \hatpi(\Hatlam)=\hatlam}
 \frac{1}{|C(h_1)|}\int_{\HatF_{\Hatlam}}\sum_{h_2\in \delta^{(2)}}
 e(\HatF_{\Hatlam})\cdot \\
 &\prod_{i=1}^ne^{2\pi\img m_i^{(h_1,h_2)}(h_1)\sigma}
 \phi(-y_i-m_i^{S^1}z+m_i^{(h_1,h_2)}(h_1)\tau-m_i^{(h_1,h_2)}(h_2),
  \tau,\sigma). 
 \end{split} 
\end{equation*}

So far we fixed a representative $h_1$ in $\delta^{(1)}$ and 
$\hatX_{\hatlam}$ in which $\delta^{(1)}$ lies. 
We now move 
$h_1$ and $\hatlam$, and sum up. Then we obtain 
the equality in Proposition
\ref{prop:fxtptformula}. 

The equality in Proposition \ref{prop:brfxtptformula} is proved 
in a parallel way. One has only to observe that 
\[ \breve{m}_i^{(h_1,h_2)}(h_1)-m_i^{(h_1,h_2)}(h_1)\equiv
  \breve{f}_{\hatlam,i}-f_{\hatlam,i} \bmod N. \]
\end{proof}

\section{Proof of main theorems; II modular property} 

For $A\in SL_2(\Z)$ we define $\hatvar^A(X;z,\tau,\sigma)$ by 
\[ \hatvar^A(X;z,\tau,\sigma)=\hatvar(X;A(z,\tau),\sigma). \] 
Similarly $\brvar^A(X;z,\tau,\sigma)$ is defined by 
\[ \brvar^A(X;z,\tau,\sigma)=\brvar(X;A(z,\tau),\sigma). \] 

\begin{lemm}\label{lemm:lift}
Let $N>1$ be an integer. 
Let $X$ be an almost complex closed orbifold 
of dimension $2n$ 
with a non-trivial $S^1$ action. Assume that there exists an orbifold line 
bundle $L$ such that $\Lambda^nTX=L^N$.  
Then the action of $S^1$ can be 
lifted to an action of some finite covering group $\tilde{S}^1$ on $L$.
\end{lemm}
\begin{proof}
$S^1$ acts on $L^N=\Lambda^nTX$. 
Locally $L$ is an $N$-fold covering of $L^N$ off 
zero-section. Hence, 
the action of $N$-fold covering $\tilde{S}^1\to S^1$ on 
$L^N$ lifts to an action on $L$. 
\end{proof}
Hereafter we assume \emph{the action of $S^1$ itself lifts 
to $L$} under the situation of Lemma \ref{lemm:lift}. 
This causes no loss of generality in view of Lemma \ref{lemm:lift}; 
we may replace $S^1$ by a suitable finite covering $\tilde{S}^1$ if 
nesessary.

\begin{lemm}
Under the situation of Lemma \ref{lemm:lift}, let $m^{S^1}$ be 
the weight of the action of $S^1$ on $L$ at $x\in F=X^{S^1}$. 
Then the tangential weights $m_i^{S^1}$ satisfy the relation 
\begin{equation}\label{eq:equilib}
 \sum_{i=1}^n m_i^{S^1}=Nm^{S^1}+l 
\end{equation} 
where $l$ belongs to $\Z$ and is independent of $x\in F$. 
\end{lemm}
\begin{proof}
The weight of 
$S^1$-action on $\Lambda^nTX$ at $x$ is $\sum_i m_i^{S^1}$. 
The weight of $S^1$-action on $L^N$ at $x$ is $Nm^{S^1}$. 
Since $\Lambda^nTX=L^N$ and the both actions 
cover the same action on $X$, 
they differ only by the fiberwise action on $\Lambda^nTX$ which is of 
the form $g\cdot v=g^lv,\ g\in S^1$, with $l\in \Z$. 
This implies \eqref{eq:equilib}. 
\end{proof}

Locally $L$ is given by a line bundle over $V_x$ with 
a lifted action of $H_x$. The group $H_x$ acts on the fiber over 
$\tilde{x}=p_x^{-1}(x)$. Let $m(h)$ be a weight of that action 
for $h\in H_x$. 
It is determined modulo integers.
Since $\Lambda^nTX=L^N$ and the the weights 
$m_i^{(h_1,h_2)}(h)$ are determined modulo $\Z$ 
we may assume that $m_i^{(h_1,h_2)}(h)$ satisfy the equality 
\begin{equation}\label{eq:weightsum}
 \sum_{i=1}^nm_i^{(h_1,h_2)}(h)=Nm(h). 
\end{equation}
If $y$ denotes 
the first Chern form of $L$, then we may also assume that
\begin{equation}\label{eq:chernsum}
 \sum_i y_i=Ny, 
\end{equation} 
since only the integral concerns in the sequel.

\begin{lemm}\label{lemm:brvarA}
Let $N>1$ be an integer. 
Let $X$ be an almost complex closed orbifold of dimension $2n$  
with a non-trivial $S^1$ action such that $|H_x|$ is relatively prime 
to $N$ for all $x\in X$. Assume that there exists an orbifold line 
bundle $L$ such that $\Lambda^nTX=L^N$. Then the modified 
orbifold elliptic genus $\brvar(X;z,\tau,\sigma)$ of level 
$N$ with $\sigma=\frac{k}{N},\ 0<k<N$, is transformed by 
$A=\begin{pmatrix} a & b \\ c & d \end{pmatrix} \in SL_2(\Z)$ in 
the following way. 
\begin{equation}\label{eq:brvarA1}
 \begin{split}
 &\quad \quad\brvar^A(X;z,\tau,\sigma)= 
 e^{\pi\img (nc(c\tau+d)\sigma^2-2clz\sigma)}\cdot \\
 &\quad \sum_{\Hatlam\in\HatLam}
 \frac{1}{|H_{\pi\circ\hatpi(\delta)}|}
 \int_{\HatF_{\Hatlam}}\sum_{(h_1,h_2)\in \delta}
 e(\HatF_{\Hatlam})e^{-2\pi\img m(h_1)dk}  
  e^{-2\pi\img(y+m^{S^1}z+m(h_2))ck}\cdot  \\
 &\prod_{i=1}^ne^{2\pi\img m_i^{(h_1,h_2)}(h_1)(c\tau+d)\sigma}
 \phi(-y_i-m_i^{S^1}z+m_i^{(h_1,h_2)}(h_1)\tau-m_i^{(h_1,h_2)}(h_2),
  \tau,(c\tau+d)\sigma),
\end{split} 
\end{equation}
where $c_1(L)=2\pi\img y$, $m^{\l h_1,h_2\r}(h_1)$ is the weight of 
$H_x$ on $L$ 
and $\delta\in \HatF_{\Hatlam}$. 
\end{lemm}
\begin{proof}
By Proposition \ref{prop:brfxtptformula} we have 
\begin{equation*}
\begin{split}
&\qquad \brvar^A(X;z,\tau,\sigma)=\sum_{\Hatlam\in\HatLam}
 \frac{1}{|H_{\pi\circ\hatpi(\delta)}|} 
 \int_{\HatF_{\Hatlam}}\sum_{(h_1,h_2)\in \delta}e(\HatF_{\Hatlam})\cdot \\
 & \prod_{i=1}^ne^{2\pi\img \breve{m}_i^{\l h_1,h_2\r}(h_1)\sigma}
 \phi(-y_i-m_i^{S^1}\frac{z}{c\tau+d}+m_i^{(h_1,h_2)}(h_1)A\tau
 -m_i^{(h_1,h_2)}(h_2),
  A\tau,\sigma) \\
&\quad =\sum_{\Hatlam\in\HatLam}\frac{1}{|H_{\pi\circ\hatpi(\delta)}|} 
 \int_{\HatF_{\Hatlam}}\sum_{(h_1,h_2)\in \delta}\frac{1}{(c\tau+d)^{r_1}}
 (\prod_{i=1}^{r_1}(c\tau+d)x_i) \cdot \\
  & \prod_{i=1}^ne^{2\pi\img \breve{m}_i^{\l h_1,h_2\r}(h_1)\sigma}
 \phi(-\frac{(c\tau+d)y_i}{c\tau+d}-m_i^{S^1}\frac{z}{c\tau+d}
 +m_i^{(h_1,h_2)}(h_1)A\tau
 -m_i^{(h_1,h_2)}(h_2),
  A\tau,\sigma).
 \end{split} 
\end{equation*}

The term of degree $2j$ of 
\[ \prod_{i=1}^{r_1}(c\tau+d)x_i \prod_{i=1}^n
 \phi(-\frac{(c\tau+d)y_i}{c\tau+d}-m_i^{S^1}\frac{z}{c\tau+d}
 +m_i^{(h_1,h_2)}(h_1)A\tau -m_i^{(h_1,h_2)}(h_2),  A\tau,\sigma) \]
is equal to $(c\tau+d)^j$ times that of
\[ \prod_{i=1}^{r_1}x_i \prod_{i=1}^n
 \phi(-\frac{y_i}{c\tau+d}-m_i^{S^1}\frac{z}{c\tau+d}
 +m_i^{(h_1,h_2)}(h_1)A\tau
 -m_i^{(h_1,h_2)}(h_2),
  A\tau,\sigma). \]
Since $\dim \HatF_{\Hatlam}=2r_1$ we obtain
\begin{equation}\label{eq:6-1}
 \begin{split}
&\qquad \brvar^A(X;z,\tau,\sigma)=\sum_{\Hatlam\in\HatLam}
 \frac{1}{|H_{\pi\circ\hatpi(\delta)}|} 
 \int_{\HatF_{\Hatlam}}\sum_{(h_1,h_2)\in \delta}e(\HatF_{\Hatlam})\cdot \\
 & \prod_{i=1}^n
 e^{2\pi\img \sum_i\breve{m}_i^{\l h_1,h_2\r}(h_1)\sigma}
 \phi(-\frac{y_i}{c\tau+d}-m_i^{S^1}\frac{z}{c\tau+d}
 +m_i^{(h_1,h_2)}(h_1)A\tau
 -m_i^{(h_1,h_2)}(h_2),
  A\tau,\sigma).
 \end{split} 
\end{equation}
Using \eqref{eq:transformula3}, \eqref{eq:equilib} and 
\eqref{eq:chernsum}, we see that
\begin{equation}\label{eq:6-2}
 \begin{split}
 &\prod_{i=1}^n
 \phi(-\frac{y_i}{c\tau+d}-m_i^{S^1}\frac{z}{c\tau+d}
 +m_i^{(h_1,h_2)}(h_1)A\tau
 -m_i^{(h_1,h_2)}(h_2),A\tau,\sigma) \\
 &=\prod_{i=1}^n
 \phi(\frac{-y_i-m_i^{S^1}z+(a\tau+b)m_i^{(h_1,h_2)}(h_1)
 -(c\tau+d)m_i^{(h_1,h_2)}(h_2)}{c\tau+d},A\tau,\sigma) \\
 &=e^{\pi\img (nc(c\tau+d)\sigma^2-2clz\sigma)}e^{-2\pi\img(y+m^{S^1}z)ck}
 e^{2\pi\img\sum_i((a\tau+b)m_i^{(h_1,h_2)}(h_1)
  -(c\tau+d)m_i^{(h_1,h_2)}(h_2))c\sigma}\cdot \\
 &\quad \prod_i \phi(-y_i-m_i^{S^1}z+(a\tau+b)m_i^{(h_1,h_2)}(h_1)
 -(c\tau+d)m_i^{(h_1,h_2)}(h_2),\tau,(c\tau+d)\sigma).
 \end{split} 
\end{equation}

In \eqref{eq:6-2} we have 
\[ \begin{split}
 ((a\tau+b)m_i^{(h_1,h_2)}(h_1)
  &-(c\tau+d)m_i^{(h_1,h_2)}(h_2))c \\
 &=-m_i^{(h_1,h_2)}(h_1)
 +(am_i^{(h_1,h_2)}(h_1)-cm_i^{(h_1,h_2)}(h_2))(c\tau+d),
 \end{split} \]
and 
\[ \begin{split}
 (a\tau+b)m_i^{(h_1,h_2)}(h_1)
 &-(c\tau+d)m_i^{(h_1,h_2)}(h_2) \\
 &=(am_i^{(h_1,h_2)}(h_1)-cm_i^{(h_1,h_2)}(h_2))\tau +
  bm_i^{(h_1,h_2)}(h_1)-dm_i^{(h_1,h_2)}(h_2).
 \end{split} \]

We consider the map $\rho:CM(H_x)\to CM(H_x)$ defined by 
\[ \rho(h_1,h_2)=(\bar{h}_1,\bar{h}_2)=(h_1^ah_2^{-c},h_1^{-b}h_2^d). \]
It is a bijection and its inverse is given by
\[ \rho^{-1}(\bar{h}_1,\bar{h}_2)=
 (\bar{h}_1^d\bar{h}_2^c,\bar{h}_1^b\bar{h}_2^a). \]
$\rho$ induces a bijection of $\hatcalC(H_x)$ onto itself which 
we shall also denote by $\rho$. It in turn induces 
a permutation $\rho$ of 
$\{\HatF_{\Hatlam}\}_{\Hatlam\in \HatLam}$ in the 
following way. If $\delta$ lies in $\HatF_{\Hatlam}$, then 
$\rho(\delta)$ lies in $\HatF_{\rho(\Hatlam)}$. 
Note that $\delta \mapsto \rho(\delta)$ defines an isomorphism of 
orbifolds from $\HatF_{\Hatlam}$ onto $\HatF_{\rho(\Hatlam)}$. 
In fact, if  
$(V_x^{t,h_1,h_2}, V_x^{t,h_1,h_2}/(C(h_1)\cap C(h_2)), C(h_1)\cap C(h_2))$ 
is a chart for $\HatF_{\Hatlam}$ and 
$(V_x^{t,\bar{h}_1,\bar{h}_2}, V_x^{t,\bar{h}_1,\bar{h}_2}/
(C(\bar{h}_1)\cap C(\bar{h}_2)), C(\bar{h}_1)\cap C(\bar{h}_2))$ 
is a chart for $\rho(\HatF_{\Hatlam})$, then 
$V_x^{t,h_1,h_2}=V_x^{t,\bar{h}_1,\bar{h}_2}$ and 
$C(h_1)\cap C(h_2)=C(\bar{h}_1)\cap C(\bar{h}_2)$. 
It follows that the identity map 
$V_x^{t,h_1,h_2}\to V_x^{t,\bar{h}_1,\bar{h}_2}$ induces 
$\rho:\HatF_{\Hatlam}\to \HatF_{\rho(\Hatlam)}$. 
$\rho:\HatF_{\Hatlam}\to \HatF_{\rho(\Hatlam)}$ respects 
$\pi\circ\hatpi$ but not necesarily 
$\hatpi$. Moreover $am_i^{(h_1,h_2)}(h_1)-cm_i^{(h_1,h_2)}(h_2)$ 
is a weight of $\bar{h}_1$ and 
$-bm_i^{(h_1,h_2)}(h_1)+dm_i^{(h_1,h_2)}(h_2)$ is a weight 
of $\bar{h}_2$. We shall use these for weights on the transformed 
sector $\rho(\HatF_{\Hatlam})$ assigned for the pair 
$(\bar{h}_1,\bar{h}_2)$ and write them 
\[ m_i^{(\bar{h}_1,\bar{h}_2)}(\bar{h}_1),\ 
  m_i^{(\bar{h}_1,\bar{h}_2)}(\bar{h}_2). \]

Thus 
\begin{equation}\label{eq:bar-1}
 \begin{split}
 ((a\tau+b)m_i^{(h_1,h_2)}(h_1)
  &-(c\tau+d)m_i^{(h_1,h_2)}(h_2))c \\
 &=-m_i^{(h_1,h_2)}(h_1)
 +m_i^{(\bar{h}_1,\bar{h}_2)}(\bar{h}_1)(c\tau+d),
 \end{split} 
\end{equation}
and 
\begin{equation}\label{eq:bar-2}
 \begin{split}
 (a\tau+b)m_i^{(h_1,h_2)}(h_1)
 &-(c\tau+d)m_i^{(h_1,h_2)}(h_2) \\
 &=m_i^{(\bar{h}_1,\bar{h}_2)}(\bar{h}_1)\tau -
   m_i^{(\bar{h}_1,\bar{h}_2)}(\bar{h}_2).
 \end{split} 
\end{equation}

Also we have
\[ m_i^{(h_1,h_2)}(h_1)=dm_i^{(\bar{h}_1,\bar{h}_2)}(\bar{h}_1)
  +cm_i^{(\bar{h}_1,\bar{h}_2)}(\bar{h}_2). \] 
Hence, by \eqref{eq:weightsum} 
\begin{equation}\label{eq:6-5}
\sum_im_i^{(h_1,h_2)}(h_1)\sigma=dkm(\bar{h}_1)
  +ckm(\bar{h}_2). 
\end{equation}

Finally, since $\sum_i\breve{m}_i^{\l h_1,h_2\r}(h_1)=N\breve{m}(h_1)$ 
by \eqref{eq:weightsum}, 
we get 
\begin{equation}\label{eq:6-4}
e^{2\pi\img \sum_i\breve{m}_i^{\l h_1,h_2\r}(h_1)\sigma}=1 
\end{equation}
for $\sigma=\frac{k}{N}$ with $0<k<N$. 

Then, by using \eqref{eq:bar-1}, \eqref{eq:bar-2} and \eqref{eq:6-5}, 
we can rewrite \eqref{eq:6-2} as follows. 
 \begin{equation}\label{eq:6-3}
 \begin{split}
 \prod_{i=1}^n
 \phi&(-\frac{y_i}{c\tau+d}-m_i^{S^1}\frac{z}{c\tau+d}
 +m_i^{(h_1,h_2)}(h_1)A\tau
 -m_i^{(h_1,h_2)}(h_2),A\tau,\sigma) \\
  =&e^{\pi\img (nc(c\tau+d)\sigma^2-2clz\sigma)}
    e^{-2\pi\img m(\bar{h}_1)dk} 
 \rho^*\left(
 e^{-2\pi\img(y+m^{S^1}z+m(\bar{h}_2))ck} \cdot
 \right. \\
 &\left. \prod_i 
 e^{2\pi\img m_i^{(\bar{h}_1,\bar{h}_2)}(\bar{h}_1)(c\tau+d)\sigma}
 \phi(-y_i-m_i^{S^1}+m_i^{\l\bar{h}_1,
 \bar{h}_2\r}(\bar{h}_1)
 -m_i^{(\bar{h}_1,\bar{h}_2)}(\bar{h}_2),\tau,(c\tau+d)\sigma)\right).
 \end{split} 
\end{equation}
We also have 
\begin{equation}\label{eq:euler}
e(\HatF_{\Hatlam})=\rho^*(e(\HatF_{\rho(\Hatlam)})).
\end{equation}

Putting \eqref{eq:6-4}, \eqref{eq:6-3} and \eqref{eq:euler} into 
\eqref{eq:6-1}, we obtain 
\begin{equation*}
 \begin{split}
 &\brvar^A(X;z,\tau,\sigma)
 = e^{\pi\img (nc(c\tau+d)\sigma^2-2clz\sigma)}\cdot \\
 &\sum_{\Hatlam\in\HatLam}
 \frac{1}{|H_{\pi\circ\hatpi(\delta)}|}
 \int_{\HatF_{\Hatlam}}\sum_{(\bar{h}_1,\bar{h}_2)\in \rho(\delta)}
 e^{-2\pi\img m(\bar{h}_1)dk} 
 \rho^*\left(e(\HatF_{\rho(\Hatlam)}) 
  e^{-2\pi\img(y+m^{S^1}z+m(\bar{h}_2))ck}
 \cdot\right.  \\
 &\left.\prod_{i=1}^ne^{2\pi\img m_i^{(\bar{h}_1,\bar{h}_2)}(\bar{h}_1)
 (c\tau+d)\sigma}
 \phi(-y_i-m_i^{S^1}z+m_i^{(\bar{h}_1,\bar{h}_2)}(\bar{h}_1)
 \tau-m_i^{(\bar{h}_1,\bar{h}_2)}(\bar{h}_2),
  \tau,(c\tau+d)\sigma)\right) \\
&=e^{\pi\img (nc(c\tau+d)\sigma^2-2clz\sigma)} \\
 &\sum_{\Hatlam\in\HatLam}
 \frac{1}{|H_{\pi\circ\hatpi(\delta)}|}
 \int_{\HatF_{\rho(\Hatlam)}}\sum_{(\bar{h}_1,\bar{h}_2)\in \rho(\delta)}
  e(\HatF_{\rho(\Hatlam)})e^{-2\pi\img m(\bar{h}_1)dk}  
  e^{-2\pi\img(y+m^{S^1}z+m(\bar{h}_2))ck}
 \cdot  \\
 &\prod_{i=1}^ne^{2\pi\img m_i^{(\bar{h}_1,\bar{h}_2)}(\bar{h}_1)
 (c\tau+d)\sigma}
 \phi(-y_i-m_i^{S^1}z+m_i^{(\bar{h}_1,\bar{h}_2)}(\bar{h}_1)
 \tau-m_i^{(\bar{h}_1,\bar{h}_2)}(\bar{h}_2),
  \tau,(c\tau+d)\sigma).
\end{split} 
\end{equation*}
Replacing $\rho(\Hatlam)$ by $\Hatlam$ and $\bar{h}_1$ and $\bar{h}_2$ 
by $h_1$ and $h_2$ respectively we obtain \eqref{eq:brvarA1}.
\end{proof}

The definition of vector bundles $\T_{\hatlam}$ and $\W_{\hatlam}$ 
depended on a parameter $\sigma$ with $\zeta=e^{2\pi\img\sigma}$. 
In case it is necessary to specify the parameter $\sigma$, we shall
write $\T_{\hatlam}(\sigma)$ and $\W_{\hatlam}(\sigma)$.  

\begin{coro}\label{coro:brvarA2}
Under the situation of Lemma \ref{lemm:brvarA} we have
\begin{equation}\label{eq:brvarA2}
 \begin{split}
 \brvar^A&(X;z,\tau,\sigma) =
 e^{\pi\img (nc(c\tau+d)\sigma^2-2clz\sigma)} 
 \sum_{\hatlam\in\hatLam} 
 e^{-2\pi\img m(h_1)dk}\cdot \\
 &e^{2\pi\img f_{\hatlam}(c\tau+d)\sigma}
 \ind\left(D_{\hatX_{\hatlam}}\otimes (L^*)^{ck}\otimes
 \T_{\hatlam}((c\tau+d)\sigma)
 \otimes\W_{\hatlam}((c\tau+d)\sigma)\right).\  
 \end{split} 
\end{equation}
\end{coro}
\begin{proof}
Apply the index formula \eqref{eq:fixpoint2} to Lemma 
\ref{lemm:brvarA}. The details of proof are similar to that of
Proposition \ref{prop:fxtptformula}. 
\end{proof}

The following Lemma is crucial for the subsequent discussion. 
\begin{lemm}\label{lemm:nopole}
The meromorphic function $\brvar^A(X;z,\tau,\sigma)$ in $z$ and $\tau$ 
has no poles at $z\in \R$. 
\end{lemm}
\begin{proof}
In view of \eqref{eq:brvarA2} it suffices to show that each
\[ \ind(D_{\hatX_{\hatlam}}\otimes (L^*)^{ck}\otimes
\T_{\hatlam}((c\tau+d)\sigma)
\otimes\W_{\hatlam}((c\tau+d)\sigma)) \]
has no poles at $z\in \R$. Furthermore 
we may replace $(L^*)^{ck}$ by an orbifold line bundle $L'$ 
and $(c\tau+d)\sigma$ by $\sigma$ without loss of generality. 
Considering 
$\ind(D_{\hatX_{\hatlam}}\otimes L'\otimes\T_{\hatlam}(\sigma)
 \otimes\W_{\hatlam}(\sigma))$ 
as a function of $z$, we write
it $\var(z)$ and expand it as a power series: 
\[ \var(z)=\sum_kb_k(z)q^{\frac{k}{r}}, \]
where $b_k(z)$ is of the form $\ind(\hat{R}'_{\hatlam,k}(\sigma))
\in R(S^1)\otimes \C$ as in \eqref{eq:var}. It follows that 
$b_k(z)$ has no poles at $z\in \R$.

On the other hand \eqref{eq:brvarA1} compared with 
\eqref{eq:brvarA2} shows that $\var(z)$ has the following expression
\begin{equation}\label{eq:varz}
\begin{split} 
 \var(z)=\zeta^{\frac{n}{2}}
 \sum_{\Hatlam\in\HatLam,
 \hat{\pi}(\HatF_{\Hatlam})\subset \hatX_{\hatlam}}
 \frac{1}{|H_{\pi\circ\hatpi(\delta)}|}& 
 \int_{\HatF_{\Hatlam}}
 \sum_{(h_1,h_2)\in \delta}e(\HatF_{\Hatlam})
 e^{2\pi\img(y+m^{S^1}z+m(h_2))}\cdot \\
 & \prod_{i=1}^n
 \phi(-y_i-m_i^{S^1}z+f_{\hatlam,i}\tau-m_i^{(h_1,h_2)}(h_2),
  \tau,\sigma),
 \end{split} 
\end{equation} 
where $y$ is the first chern form of $L'$, and 
$m^{S^1}$ and $m(h)$ are the weights of the actions of 
$S^1$ and $h\in H_x$ on $L'_x$ respectively. 
If $\sum_kb_{k,\Hatlam}(z)q^{\frac{k}{r}}$ is the contribution 
from $\Hatlam\in\HatLam$ to \eqref{eq:varz}, then  
\[ b_k(z)= \sum_{\Hatlam\in\HatLam}b_{k,\Hatlam}(z). \]
The term of degree zero in the denominator of 
$\phi(-y_i-m_i^{S^1}z+f_{\hatlam,i}\tau-m_i^{(h_1,h_2)}(h_2),
  \tau,\sigma)$ is 
\[ \begin{split}(1-&e^{2\pi\img(-m_i^{S^1}z-m_i^{(h_1,h_2)}(h_2))}
 q^{f_{\hatlam,i}\tau})\cdot \\
 &\prod_{k=1}^\infty
 (1-e^{2\pi\img(-m_i^{S^1}z-m_i^{(h_1,h_2)}(h_2))}
 q^{f_{\hatlam,i}+k})
 (1-e^{2\pi\img(m_i^{S^1}z+m_i^{(h_1,h_2)}(h_2))}
 q^{-f_{\hatlam,i}+k}). 
 \end{split} \]
Therefore poles of $b_{k,\Hatlam}(z)$ lie in $\R$ but they 
cancel out in the sum $b_k(z)=\sum_{\Hatlam\in\HatLam}b_{k,\Hatlam}(z)$ 
as noted above. Assume $z_0\in \R$ is a pole of $\var(z)$. 
Then there is an open set $U$ containing 
$z_0$ such that 
the power series $\sum_kb_{k,\Hatlam}(z)q^{\frac{k}{r}}$ converges 
uniformly on any compact set in $U\setminus \{z_0\}$ 
and $b_k(z)=\sum_{\Hatlam\in\HatLam}b_{k,\Hatlam}(z)$ 
is holomorphic in $U$. 
In such a situation one can conclude that 
$\var(z)=\sum_kb_k(z)q^{\frac{k}{r}}$ has no poles in $\R$. We 
refer to Lemma in Section 7 of \cite{Hi}. See also Section 5 of 
\cite{HM}.
\end{proof}

We now proceed to the proof of Theorem \ref{thm:brvar}. We regard 
$\brvar(X;z,\tau,\sigma)$ as a meromorphic function of $z$. 
By the transformation law \eqref{eq:transformula2} $\phi(z,\tau,\sigma)$ 
is an elliptic function in $z$ with respect to the lattice 
$\Z\cdot N\tau\oplus \Z$ for $\sigma=\frac{k}{N}$ with $0<k<N$. 
Hence the equivariant modified orbifold elliptic genus 
$\brvar(X;z,\tau,\sigma)$ of level $N$ is also an elliptic function 
in $z$. Thus, in order to show that $\brvar(X;z,\tau,\sigma)$ is a 
constant it suffices to show that it does not have poles. 

Assume that $z$ 
is a pole. Then $1-t^mq^r\alpha =0$ for some integer $m\not=0$, some
rational number $r$ and a root of unity $\alpha$. 
Consequently there are intergers $m_1\not=0$ and $k_1$ such that
$ m_1z+k_1\tau \in\Z$. Then there is an element 
$A= \begin{pmatrix} a & b \\ c & d \end{pmatrix} \in SL_2(\Z)$ such 
that
\[ \frac{z}{c\tau+d}\in \R. \]
Since
\[ \brvar(X;z,\tau,\sigma)=\brvar(X;A^{-1}(\frac{z}{c\tau+d},A\tau),\sigma)=
 \brvar^{A^{-1}}(X;\frac{z}{c\tau+d},A\tau,\sigma), \]
$\brvar^{A^{-1}}(X;\frac{z}{c\tau+d},A\tau,\sigma)$ must have a pole 
at $\frac{z}{c\tau+d}\in \R$. But this contradicts with Lemma 
\ref{lemm:nopole}. This contradiction proves that $\brvar(X;z,\tau,\sigma)$ 
can not have a pole and Theorem \ref{thm:brvar} follows.

\begin{prop}\label{prop:vanishing}
Let $X$ be an almost complex closed orbifold 
of dimension $2n$ 
with a non-trivial $S^1$ action. Let $N>1$ be an integer 
relatively prime to the orders of all isotropy groups $|H_\lam|$. 
Assume that there exists an orbifold line 
bundle $L$ such that $\Lambda^nTX=L^N$. If the number $l$ in 
\eqref{eq:equilib} is relatively prime to $N$, then the modified orbifold 
elliptic genus $\brvar(X)$ of level $N$ vanishes. 
\end{prop}
\begin{proof}
By Theorem \ref{thm:brvar} the equivariant modified elliptic genus 
$\brvar_t(X;\tau,\sigma))=\brvar(X;z,\tau,\sigma)$ is constnat and 
equal to $\brvar(X;\tau,\sigma)$. By 
\eqref{eq:transformula2}, Proposition \ref{prop:brfxtptformula} and 
\eqref{eq:equilib} we have 
\[ \brvar(X;\tau,\sigma)=\brvar(X;z+\tau,\tau,\sigma)=
 \zeta^l\brvar(X;z,\tau,\sigma)=
 \zeta^l\brvar(X;\tau,\sigma). \]
Since $l$ is relatively prime to $N$, $\zeta^l$ is not
equal to $1$. Hence $\brvar(X;\tau,\sigma)$ must vanish.
\end{proof}

By similar calculations to the ones used in the proof of Lemma 
\ref{lemm:brvarA} we obtain the following 

\begin{lemm}\label{lemm:hatvarA}
Let $X$ be an almost complex closed orbifold 
of dimension $2n$ with a non-trivial $S^1$ action and 
let $N>1$ be an integer. Assume that there exists a genuine line 
bundle $L$ such that $\Lambda^nTX=L^N$. Then the 
orbifold elliptic genus $\hatvar(X;z,\tau,\sigma)$ of level 
$N$ with $\sigma=\frac{k}{N},\ 0<k<N$, is transformed by 
$A=\begin{pmatrix} a & b \\ c & d \end{pmatrix} \in SL_2(\Z)$ in 
the following way. 
\begin{equation}\label{eq:hatvarA2}
 \begin{split}
 &\hatvar^A(X;z,\tau,\sigma) = 
 e^{\pi\img (nc(c\tau+d)\sigma^2-2clz\sigma)}\sum_{\Hatlam\in\HatLam}
 \frac{1}{|H_{\pi\circ\hatpi(\delta)}|}
 \int_{\HatF_{\Hatlam}}\sum_{(h_1,h_2)\in \delta}e(\HatF_{\Hatlam}) 
  e^{-2\pi\img(y+m^{S^1}z)ck}\cdot \\
 &\prod_{i=1}^ne^{2\pi\img m_i^{(h_1,h_2)}(h_1)(c\tau+d)\sigma}
 \phi(-y_i-m_i^{S^1}z+m_i^{(h_1,h_2)}(h_1)\tau-m_i^{(h_1,h_2)}(h_2),
  \tau,(c\tau+d)\sigma).
\end{split} 
\end{equation}
\end{lemm}
\begin{note}
Under the assumption of Lemma \ref{lemm:hatvarA} (including the 
case $N=1$) each $f_{\hatlam}$ 
is an integer. In fact $\sum_im_i^{(h_1,h_2)}(h)=Nm(h)$ is an integer 
as the weight of $h$ on a genuine line bundle $L^N$ by \eqref{eq:weightsum}. 
$f_{\hatlam}$ is congruent to $\sum_im_i^{(h_1,h_2)}(h)$ for 
$[h]\in \hatX_{\hatlam}$. 
\end{note}
The proof of Theorem \ref{thm:hatvar} can be given in a similar 
way to that of Theorem \ref{thm:brvar} by using Lemma 
\ref{lemm:hatvarA}. In fact, under the assumption that the 
orbifold line bundle $L$ is genuine, the term 
$e^{-2\pi\img(y+m^{S^1}z)ck}$ in the integrand of 
\eqref{eq:hatvarA2} is equal to 
$e^{-2\pi\img(y+m^{S^1}z+m(h_2))ck}$ 
so that we have 
\begin{equation*}\label{eq:hatvarA3}
 \begin{split}
 \hatvar^A&(X;z,\tau,\sigma) =
 e^{\pi\img (nc(c\tau+d)\sigma^2-2clz\sigma)}\cdot \\
 & \sum_{\hatlam\in\hatLam} 
 e^{2\pi\img f_{\hatlam}(c\tau+d)\sigma}
 \ind\left(D_{\hatX_{\hatlam}}\otimes (L^*)^{ck}\otimes\T_{\hatlam}
 ((c\tau+d)\sigma)
 \otimes\W_{\hatlam}((c\tau+d)\sigma)\right).  
 \end{split} 
\end{equation*} 
The rest of the proof is entirely similar to that of Theorem 
\ref{thm:brvar}. One sees that $\hatvar(X;z,\tau,\sigma)$ can 
not have a pole and consequently it is constant. 

\begin{prop}\label{prop:hatvanishing}
Let $X$ be an almost complex closed orbifold 
of dimension $2n$ 
with a non-trivial $S^1$ action and let $N>1$ be an integer. 
Assume that there exists a genuine line 
bundle $L$ such that $\Lambda^nTX=L^N$. If the number $l$ in 
\eqref{eq:equilib} is relatively prime to $N$, then the orbifold 
elliptic genus $\hatvar(X)$ of level $N$ vanishes. 
\end{prop}
\begin{proof}
In a similar way to the proof of Proposition \ref{prop:vanishing}, 
we have  
\[ \hatvar(X;\tau,\sigma)=\hatvar(X;z+\tau,\tau,\sigma)=
  \zeta^l\hatvar(X;z,\tau,\sigma)=
  \zeta^l\hatvar(X;\tau,\sigma). \]
Since $\zeta^l\not=1$, $\hatvar(X;\tau,\sigma)$ must vanish. 
\end{proof}

The proof of Theorem \ref{thm:hatvartorsion} goes as follows. 
Assume that $\Lambda^nTX$ is trivial as an 
orbifold line bundle. Then we can take Chern form of 
$\Lambda^nTX$ to be zero form. Hence we can assume that 
$\sum_{i=1}^ny_i=0$. We also see that 
\begin{equation}\label{eq:equilib2}
 \sum_{i=1}^nm_i^{S^1}=l  
\end{equation}
at each point $x\in F$, where $l$ belongs to $\Z$ and 
is independent of $x$. We may also assume that 
$\sum_im_i^{(h_1,h_2)}(h)=0$ as in \eqref{eq:weightsum}. 
Then, a similar calculation 
to \eqref{eq:hatvarA2} yields  
\begin{equation*}\label{eq:hatvarA}
 \begin{split}
 &\hatvar^A(X;z,\tau,\sigma) = 
 e^{\pi\img (nc(c\tau+d)\sigma^2-2clz\sigma)}\sum_{\Hatlam\in\HatLam}
 \frac{1}{|H_{\pi\circ\hatpi(\delta)}|}
 \int_{\HatF_{\Hatlam}}\sum_{(h_1,h_2)\in \delta}e(\HatF_{\Hatlam}) \\
 &\prod_{i=1}^ne^{2\pi\img m_i^{(h_1,h_2)}(h_1)(c\tau+d)\sigma}
 \phi(-y_i-m_i^{S^1}z+m_i^{(h_1,h_2)}(h_1)\tau-m_i^{(h_1,h_2)}(h_2),
  \tau,(c\tau+d)\sigma).
\end{split} 
\end{equation*}
Then, as in Corollary \ref{coro:brvarA2}, we have 
\begin{equation*}
 \begin{split}
 &\hatvar^A(X;z,\tau,\sigma) = \\
 & e^{\pi\img (nc(c\tau+d)\sigma^2-2clz\sigma)} 
 \sum_{\hatlam\in\hatLam} 
 e^{2\pi\img f_{\hatlam}(c\tau+d)\sigma}
 \ind\left(D_{X_{\hatlam}}\otimes\T_{\hatlam}((c\tau+d)\sigma)
 \otimes\W_{\hatlam}((c\tau+d)\sigma)\right).
\end{split} 
\end{equation*} 
Let $N$ be an integer greater than $1$. Then
an entirely similar argument to the one in the proof of Theorem 
\ref{thm:brvar} proves that $\hatvar(X;z,\tau,\sigma)$ is constant 
for $\sigma=\frac{k}{N}, 0<k<N$. Since this is true for any integer $N>1$ 
and $\sigma=\frac{k}{N}$, 
$\hatvar(X;z,\tau,\sigma)$ must be constant. 

\begin{prop}\label{prop:trivialcone}
Let $X$ be an almost complex closed orbifold 
of dimension $2n$ 
with a non-trivial $S^1$ action. Assume that $\Lambda^nTX$ is  
trival. If the number $l$ in 
\eqref{eq:equilib2} is not zero, then  $\hatvar(X)$ vanishes. 
\end{prop}
\begin{proof}

If $l$ is not equal to $0$, then $\zeta^l\not=1$ for 
$\sigma=\frac{k}{N}, 0<k<N$, with $N$ relatively prime to $l$. 
Then a similar argument to the proof of 
Proposition \ref{prop:vanishing} proves that $\hatvar(X;\tau,\sigma)=0$ 
for such $\sigma$. This implies that $\hatvar(X;\tau,\sigma)$ must vanish. 
\end{proof}

\section{Orbifold $T_y$-genus}

Let $X$ be an almost complex closed orbifold. The $T_y$-genus 
$T_y(X)$ is defined to be the index
\begin{equation*}\label{eq:Ty}
 T_y(X)= \ind (D\otimes \Lambda_yT^*X)\in \Z[y,y^{-1}].
\end{equation*}
We further consider the orbifold $T_y$-genus
\begin{equation}\label{eq:hatTy}
 \hatT_y(X)=\sum_{\hatlam\in\hatLam}(-y)^{f_{\hatlam}}T_y(\hatX_{\hatlam})= 
 \sum_{\hatlam\in\hatLam}(-y)^{f_{\hatlam}} 
 \ind (D_{\hatX_{\hatlam}}\otimes \Lambda_yT^*\hatX_{\hatlam}).
\end{equation}
In case the orders of all isotropy groups $H_x$ are relatively 
prime to an integer $N>1$ the modified orbifold $T_y$-genus 
$\brT_y(X)$ of level $N$ is defined by 
\begin{equation}\label{eq:brTy}
 \brT_y(X)=\sum_{\hatlam\in\hatLam}(-y)^{\brf_{\hatlam}}T_y(\hatX_{\hatlam})
 =\sum_{\hatlam\in\hatLam}(-y)^{\brf_{\hatlam}} 
 \ind (D_{\hatX_{\hatlam}}\otimes \Lambda_yT^*\hatX_{\hatlam}), \ 
\end{equation}
where $-y=e^{2\pi\img\frac{k}{N}},\ 0<k<N$.
\begin{note}
$T_y(X)$ and $\hatT_y(X)$ are the constant terms in the 
power series expansions \eqref{eq:var} of the elliptic genus $\var(X)$ 
and orbifold elliptic genus $\hatvar(X)$ with $\zeta$ replaced 
by $-y$. Similarly $\brT_y(X)$ is the 
constant term of the $q$-expansion of $\brvar(X)$ with $\zeta$ replaced 
by $-y$ when $\zeta=e^{2\pi\img\frac{k}{N}},\ 0<k<N$. 
\end{note}

When a compact group $G$ acts on $X$ one can consider the corresponding 
equivariant genera. 
It is known \cite{Ko} that the $T_y$-genus is rigid for closed manifolds with 
compact connected group action. 
\begin{prop}\label{T_y}
The $T_y$-genus, the $\hatT_y$-genus and the $\brT_y$-genus are rigid for 
closed orbifolds with compact connected Lie group action. 
\end{prop}
\begin{proof}
It is enough to prove the statement for $T_y$-genus in view of 
\eqref{eq:hatTy} and \eqref{eq:brTy}. We may further assume that 
$G$ is the circle group $S^1$. Let $T_{y,t}(X)\in \Z[y,y^{-1}]\otimes
R(S^1)$ be the equivariant 
$T_y$-genus of $X$. We shall use the notations in Section 5 and 
put $\zeta=-y$. By the fixed point formula 
\eqref{eq:fixpoint2} we have 
\[ \begin{split}
 &T_{y,t}(X)= \\ 
 &\sum_{\hatlam\in \hatLam_F}
 \frac{1}{|H_{\pi(\hatlam)}|}\int_{\hatF_{\hatlam}}\sum_{h\in \gamma_x}
 Td(\hatF_{\hatlam})\prod_{i=1}^{r_1}(1-\zeta e^{-2\pi\img y_i})
 \prod_{i=r_1+1}^{n}
 \frac{(1-\zeta t^{-m_i^{S^1}}e^{-(2\pi\img y_i+m_i^{(1,h)}(h))})}
 {(1-t^{-m_i^{S^1}}e^{-(2\pi\img y_i+m_i^{(1,h)}(h))})}. 
 \end{split} \]
Note that $r_1$ and $m_i^{S^1}$ depend on $\hatlam$, and $m_i^{S^1}\not=0$ for 
$r_1+1\leq i\leq n$. Put 
\[ \mu(\hatlam)=\#\{i\mid m_i^{S^1}>0\}. \] 
Then $\#\{i\mid m_i^{S^1}<0\}=n-r_1-\mu(\hatlam)$. 
We regard 
\[ \alpha_{\hatlam}(t)=\frac{1}{|H_{\pi(\hatlam)}|}\int_{\hatF_{\hatlam}}
 \sum_{h\in \gamma_x}
 Td(\hatF_{\hatlam})\prod_{i=1}^{r_1}(1-\zeta e^{-2\pi\img y_i})
 \prod_{i=r_1+1}^{n}
 \frac{(1-\zeta t^{-m_i^{S^1}}e^{-(2\pi\img y_i+m_i^{(1,h)}(h))})}
 {(1-t^{-m_i^{S^1}}e^{-(2\pi\img y_i+m_i^{(1,h)}(h))})} \]
as a rational function of $t$.
Then we see easily that
\begin{align*} 
 \alpha_{\hatlam}(0)&=\frac{1}{|H_{\pi(\hatlam)}|}\int_{\hatF_{\hatlam}}
 \sum_{h\in \gamma_x}
 Td(\hatF_{\hatlam})\prod_{i=1}^{r_1}(1-\zeta e^{-2\pi\img y_i})
 \zeta^{\mu(\hatlam)}, \\
 \alpha_{\hatlam}(\infty)&=\frac{1}{|H_{\pi(\hatlam)}|}\int_{\hatF_{\hatlam}}
 \sum_{h\in \gamma_x}
 Td(\hatF_{\hatlam})\prod_{i=1}^{r_1}(1-\zeta e^{-2\pi\img y_i})
 \zeta^{n-r_1-\mu(\hatlam)}.
\end{align*}
Thus $T_{y,t}(X)=\sum_{\hatlam}\alpha_{\hatlam}(t)$ takes finite values 
at $t=0$ and $t=\infty$. Since $T_{y,t}(X)$ belongs to 
$\Z[y,y^{-1}]\otimes R(S^1)=\Z[y,y^{-1}]\otimes\Z[t,t^{-1}]$, it 
must be a constant which is equal to $T_y(X)$. 
\end{proof}
In the above proof we have in fact proved the following
\begin{coro}
\[ \begin{split}
T_y(X)&=\sum_{k=0}^{n-r_1}(-y)^k\sum_{\hatlam:\mu(\hatlam)=k}
 T_y(\hatF_{\hatlam}) \\
 &=\sum_{k=0}^{n-r_1}(-y)^k\sum_{\hatlam:\mu(\hatlam)=n-r_1-k}
 T_y(\hatF_{\hatlam})
\end{split} \]
\end{coro}

Note that $T_y(X)$ is a polynomial 
in $-y$ of degree at most $n=\frac{\dim X}{2}$ with integer coefficients. 
\begin{lemm}\label{lemm:atmostn}
If $\Lambda^nTX$ is a genuine line bundle, 
then $\hatT_y(X)$ is 
a polynomial in $-y$ of degree at most n with integer coefficients. 
\end{lemm}
\begin{proof}
Each $f_{\hatlam}$ is an integer by Note after Lemma
\ref{lemm:hatvarA}. Therefore $\hatT_y(X)$ is a polynomial in $-y$ 
with integer coefficients. If $\dim \hatX_{\hatlam}=2k<2n$, then 
$f_{\hatlam}=\sum_{i=1}^{n-k}f_{\hatlam,i}$ with 
$0<f_{\hatlam,i}<1$. Therefore $0<f_{\hatlam}<n-k$ and the 
degree of $(-y)^{f_{\hatlam}}T_y(\hatX_{\hatlam})$ is less than 
$n-k+k=n$. 
\end{proof}

\begin{prop}\label{prop:hatT}
Let $X$ be an almost complex closed orbifold 
of dimension $2n$ 
with a non-trivial action of a compact connected Lie group $G$. 
Let $N>1$ be an integer. 
Assume that there exists a genuine line 
bundle $L$ with a lifted action of $G$ such that 
$\Lambda^nTX=L^N$. If the number $l$ in 
\eqref{eq:equilib} is relatively prime to $N$, then the orbifold 
$\hatT_y$-genus $\hatT_y(X)$ is a polynomial in $-y$ divisible by
\[ \sum_{k=0}^{N-1}(-y)^k .\]
Moreover if $\hatT_y(X)\not=0$, then $N\leq n+1$. 
\end{prop}
\begin{proof}
By Proposition \ref{prop:hatvanishing} the orbifold elliptic 
genus $\hatvar(X)$ of level $N$ vanishes. In particular its degree $0$
term $\hatT_y(X)=\sum_{\hatlam}(-y)^{f_{\hatlam}}T_y(\hatX_{\hatlam})$ 
vanishes for $-y=e^{2\pi\img\frac{k}{N}},\ 0<k<N$. 
Since $\hatT_y(X)$ is a polynomial with integer coefficients in $-y$ 
it must be divisible by $\sum_{k=0}^{N-1}(-y)^k$. 

Assume that $\hatT_y(X)\not=0$. Since its degree is at most $n$ 
by Lemma \ref{lemm:atmostn} and 
it is divisible by $\sum_{k=0}^{N-1}(-y)^k$, we must have $N-1\leq n$. 
\end{proof}

\begin{rem}
Suppose that the multiplicity of $X$ is equal to $1$ 
and $\Lambda^nTX$ is a genuine line bundle. 
The constant term of $\hatT_y(X)$ considered as 
a polynomial of $-y$ is equal to the constant term $T_0(X)$ of 
$T_y(X)$ which is nothing but the Todd genus of $X$, 
since $f_{\hatlam}>0$ for twisted sectors. Thus, if $T_0(X)\not=0$, 
then $\hatT_y(X)$ does not vanish.
\end{rem}

Situations like Proposition \ref{prop:hatT} occur 
when an $n$-dimensional torus acts 
on $X$, cf. \cite{Ha}. 

Finally as corollaries of Proposition \ref{prop:vanishing} 
and Proposition \ref{prop:trivialcone} 
we have 
\begin{prop}
Let $X$ be an almost complex closed orbifold 
of dimension $2n$ 
with a non-trivial action of a compact connected Lie group $G$. 
Let $N>1$ be an integer 
relatively prime to the orders of all isotropy groups $|H_\lam|$. 
Assume that there exists an orbifold line 
bundle $L$ wtih a lifted action of $G$ such that 
$\Lambda^nTX=L^N$. If the number $l$ in 
\eqref{eq:equilib} is relatively prime to $N$, then the modified 
$\brT_y$-genus $\brT_y(X)$ of level $N$ vanishies for 
$-y=e^{2\pi\img\frac{k}{N}},\ 0<k<N$. 
\end{prop}

\begin{prop}
Let $X$ be an almost complex closed orbifold 
of dimension $2n$ 
with a non-trivial action of a compact connected Lie group $G$. 
Assume that $\Lambda^nTX$ is trivial as an orbifold line bundle. 
If the number $l$ in 
\eqref{eq:equilib2} is not equal to $0$, then the orbifold 
$\hatT_y$-genus $\hatT_y(X)$ vanishes. 
\end{prop}

\section{Appendix:stably almost complex orbifolds}
\label{sec:app}

A stably almost complex orbifold is an oriented orbifold with a 
stably almost complex structure on the tangent bundle. 
More precisely, let $\U$ be an orbifold atlas of $X$. Let $k$ 
be a positive integer and let $\bm{k}$ denote 
the trivial real vector bundle of dimension $k$ endowed 
with the standard orientation. Then 
a complex structure on $TV\oplus \bm{k}$ is given for each 
chart $(V,U,H)\in \U$ in such a way 
that it is compatible with 
the orientation of $TV$ followed by that of $\bm{k}$, and it is
preserved by each element $h\in H$. 
Here $h$ acts trivially on $\bm{l}$. It is also required that 
theses complex structures are compatible with injections of charts. 
Two complex structures
$TV\oplus \bm{k}$ and $TV\oplus \bm{k'}$ are considered as 
equivalent if the Whitney sums of $\bm{l}_\C$ and 
$\bm{l'}_\C$ to them give the same complex
structures on $TV\oplus (\bm{k+2l})=TV\oplus (\bm{k'+2l'})$, where 
$\bm{l}_\C$ and $\bm{l'}_\C$ are trivial 
complex vector bundles of dimension $l$ and $l'$ such that
$k+2l=k'+2l'$. An equivalence class is the stably almost 
complex structure. 
Note that each sector $\hatX_{\hatlam}$ is also a stably 
almost complex orbifold. 
In fact, for $\gamma_x=[h]\in \hatX_{\hatlam}$ with 
$h\in H_x\subset H$, the complex vector bundle 
$(TV\oplus\bm{k})^h=TV^h\oplus \bm{k}$  gives the stably 
almost complex structure on $\hatX_{\hatlam}$. 
The normal bundle 
$W_{\hatlam}$ of the immersion $\pi:\hatX_{\hatlam}\to X$ has 
a canonical complex structure and the eigen-bundle 
decomposition by the action of $h$ just 
like \eqref{eq:eigendecomp2}. 

In general let $X$ be an orbifold. For an orbifold complex vector 
bundle $W$ of rank $l$ we put $\tilde{W}=W-\bm{l}_{\C}\in K_{orb}(X)$. 
If $X$ is a stably almost complex orbifold, 
$\widetilde{TX\oplus \bm{k}}$ is a well-defined element in $K_{orb}(X)$ 
which we simply denote by $\tilde{T}X$. Similarly we put 
$\tilde{T}^*X=\widetilde{T^*X\oplus\bm{k}}$. 

We then define
formal vector bundles $\tilde{\T}=\tilde{\T}(\sigma)$ and 
$\tilde{\T}_{\hatlam}=\tilde{\T}_{\hatlam}(\sigma)$
 by
\begin{align*}
 \tilde{\T}(\sigma)=&\Lambda_{-\zeta}\tilde{T}^*X\otimes
  \bigotimes_{k=1}^\infty\left(
   \Lambda_{-\zeta q^k}\tilde{T}^*X\otimes 
   \Lambda_{-\zeta^{-1}q^k}\tilde{T}X
   \otimes S_{q^k}\tilde{T}^*X\otimes S_{q^k}\tilde{T}X \right), \\
 \tilde{\T}_{\hatlam}(\sigma)=&
   \Lambda_{-\zeta}\tilde{T}^*\hatX_{\hatlam}\otimes
   \bigotimes_{k=1}^\infty\left(
   \Lambda_{-\zeta q^k}\tilde{T}^*\hatX_{\hatlam}\otimes 
   \Lambda_{-\zeta^{-1}q^k}\tilde{T}\hatX_{\hatlam}
   \otimes S_{q^k}\tilde{T}^*\hatX_{\hatlam}\otimes 
    S_{q^k}\tilde{T}\hatX_{\hatlam} \right).
\end{align*}
It is to be noted that 
\[ \Lambda_\lambda \tilde{W}=\Lambda_\lambda W/(1+\lambda)^{\rank W}, 
\quad S_\lambda\tilde{W}=S_\lambda W/(1+\lambda)^{\rank W}, \]
because of multiplicative 
properties of total exterior power operation and symmetric 
power operation. It follows that
\[ 
 \tilde{\T}(\sigma)=\T(\sigma)/(-\zeta^{\frac12}
 \Phi(\sigma,\tau))^{\frac{\dim X}{2}}, \quad 
 \tilde{\T}_{\hatlam}(\sigma)=\T(\sigma)_{\hatlam}
 /(-\zeta^{\frac12}\Phi(\sigma,\tau))^{\frac{\dim X_{\hatlam}}{2}}, \]  
when $X$ is an almost complex orbifold. 
Since we have eigen-bundle decomposition \eqref{eq:eigendecomp2} 
as in the case of almost complex orbifold, we can define 
$\W_{\hatlam,i}=\W_{\hatlam,i}(\sigma)$ and 
$\W_{\hatlam}=\W_{\hatlam}(\sigma)$ just as in section 3. We 
then define 
\begin{equation*}
  \tilde{\W}_{\hatlam}=\W_{\hatlam}
  /(-\zeta^{\frac12}\Phi(\sigma,\tau))^{\rank W_{\hatlam}}. 
\end{equation*} 

On a stably almost complex orbifold a spin-c Dirac operator can 
be introduced as in the case of almost complex manifolds. 
It is an operator of the same form as \eqref{eq:Dirac} with
\[ E^+=\bigoplus_{i:even}\Lambda^i(TX\oplus \bm{k}) \ \text{and}\ 
   E^-=\bigoplus_{i:odd}\Lambda^i(TX\oplus \bm{k}). \]

Hereafter we assume that $X$ is a stably almost complex orbifold 
of dimension $2n$, and 
a complex vector bundle structure of rank $n+s$ 
is given on $T'X=TX\oplus \bm{2s}$ for some $s$. 

We define the stabilized elliptic genus $\var_{st}(X)$ and 
stabilized orbifold elliptic genus $\hatvar_{st}(X)$ of $X$ 
by 
\begin{equation*}
 \begin{split}
 \var_{st}(X)=&(-1)^n\ind(D\otimes \tilde{\T}(\sigma)) \\
 \hatvar_{st}(X)=&(-1)^n\sum_{{\hatlam}\in \hatLam}
 \zeta^{f_{\hatlam}}\ind(D_{\hatX_{\hatlam}}
 \otimes\tilde{\T}_{\hatlam}(\sigma)
  \otimes \tilde{\W}_{\hatlam}(\sigma)). 
\end{split} 
\end{equation*}
When $X$ is an almost complex orbifold of dimension $2n$ we have 
\begin{equation}\label{eq:stable}
 \var(X)=\var_{st}(X)\Phi(\sigma,\tau)^n,\quad 
 \hatvar(X)=\hatvar_{st}(X)\Phi(\sigma,\tau)^n.
\end{equation}
We may define $\var(X)$ and $\hatvar(X)$ for stably almost 
complex orbifold $X$ by \eqref{eq:stable}. 

Suppose that $N>1$ is an integer relatively prime to 
every $|H_x|$. We define modified stabilized orbifold elliptic 
genus $\brvar_{st}(X)$ by 
\[ \brvar_{st}(X)=(-1)^n\sum_{{\hatlam}\in \hatLam}
 \zeta^{\brf_{\hatlam}}\ind(D_{\hatX_{\hatlam}}
 \otimes\tilde{\T}_{\hatlam}(\sigma)
  \otimes \tilde{\W}_{\hatlam}(\sigma)), \]
and set $\brvar(X)=\brvar_{st}(X)\Phi(\sigma,\tau)^n$. 

The cohomology classes of $c(T'X)\in H^*(X,\Z_X)$ 
and $Td(T'X)\in H^*(X,\Q_X)$ are well-defined classes 
depending only on $\tilde{T}X$, where
\[ Td(T'X)=\det(\frac{\Gamma(T'X)}{1-e^{-\Gamma(T'X)}})=
 \prod_{i=1}^{n+s}\frac{x_i}{1-e^{-x_i}}, \]
with $c(T'X)=\prod_i(1+x_i)$ as before. 

When a compact connected Lie group $G$ acts on $X$ it is always 
assumed that the action preserves the stably almost complex 
structure. 
Then we can consider equivariant genera corresponding 
$\hatvar(X)$ etc.
The fixed point set $X^G$ is a stably almost complex 
orbifold. Vergne's fixed point formulas \eqref{eq:fixpoint} and 
\eqref{eq:fixpoint2} still hold for stably almost complex orbifolds 
by replacing $Td(\hatF_{\hatlam})$ by $Td(T'\hatF_{\hatlam})$. 

Suppose that $N>1$ is an integer and there is an orbifold line 
bundle $L$ such that $\Lambda^{n+s}T'X=L^N$. 
This condition is equivalent to saying that the first orbifold 
Chern class $c_1(\tilde{T}X)\in H^2(X,\Z_X)$ is divisible by $N$. 
Similarly the condition that $\Lambda^{n+s}T'X$ is trivial 
means that $c_1(\tilde{T}X)\in H^2(X,\Z_X)$ vanishes. 
Theorems \ref{thm:brvar}, \ref{thm:hatvar} and \ref{thm:hatvartorsion} 
have meanings for stably almost complex manifolds by replacing 
$\Lambda^nTX$ by $\Lambda^{n+s}T'X$, 
and they in fact hold 
in this extended sense. Similarly Propositions \ref{prop:vanishing}, 
\ref{prop:hatvanishing} and \ref{prop:trivialcone} hold for stably 
almost complex orbifolds. 

Proofs are almost verbatim. 
We work with stabilized genera 
$\hatvar_{st}(X)$ and so on. We put 
\[ \phi_{st}(z,\tau,\sigma)= \phi(z,\tau,\sigma)/\Phi(\sigma,\tau)
 =\Phi(z+\sigma,\tau)/\Phi(z,\tau)\Phi(\sigma,\tau). \]
$\phi_{st}(z,\tau,\sigma)$ satisfies the following transformation
law. 
\begin{equation*}\label{eq:sttransformula3}
 \begin{split}
 &\phi_{st}(A(z,\tau),\sigma)=(c\tau+d)e^{2\pi\img cz\sigma}
 \phi_{st}(z,\tau,(c\tau+d)\sigma) \ \ \text{for}\ A=
 \begin{pmatrix} a & b \\ c & d \end{pmatrix}\in SL_2(\Z), \\
 &\phi_{st}(z+m\tau+n,\tau,\sigma)=
 e^{-2\pi\img m\sigma}\phi_{st}(z,\tau,\sigma)\ \ 
 \text{for} \ m,n\in \Z. 
 \end{split} 
\end{equation*}
 
We put $e(T'\HatF_{\Hatlam})=\prod_{i=1}^{r_1}x_i$ where 
$T'\HatF_{\Hatlam}=T\HatF_{\Hatlam}\oplus\bm{2s}$ and 
$r_1=\rank_{\C}(T'\HatF_{\Hatlam})=
\frac{\dim\HatF_{\Hatlam}}{2}+s$. We use the same notions and 
conventions as in \eqref{eq:split2}, \eqref{eq:weights} and 
\eqref{eq:weights2}. 
Then \eqref{eq:fxtptformula} and \eqref{eq:brfxtptformula} 
are replaced by 
\begin{equation}\label{eq:stfxtptformula}
\begin{split} 
 \hatvar_{st}&(X;z,\tau,\sigma)=\sum_{\Hatlam\in\HatLam}
 \frac{1}{|H_{\pi\circ\hatpi(\delta)}|} 
 \int_{\HatF_{\Hatlam}}\sum_{(h_1,h_2)\in \delta}
 e(T'\HatF_{\Hatlam})\cdot \\
 & \prod_{i=1}^ne^{2\pi\img m_i^{(h_1,h_2)}(h_1)\sigma}
 \phi_{st}(-y_i-m_i^{S^1}z+m_i^{(h_1,h_2)}(h_1)\tau-m_i^{(h_1,h_2)}(h_2),
  \tau,\sigma)
 \end{split} 
\end{equation}
and
\begin{equation}\label{eq:stbrfxtptformula}
\begin{split} 
 \brvar_{st}&(X;z,\tau,\sigma)=\sum_{\Hatlam\in\HatLam}
 \frac{1}{|H_{\pi\circ\hatpi(\delta)}|} 
 \int_{\HatF_{\Hatlam}}\sum_{(h_1,h_2)\in \delta}
 e(T'\HatF_{\Hatlam})\cdot \\
 & \prod_{i=1}^ne^{2\pi\img \breve{m}_i^{\l h_1,h_2\r}(h_1)\sigma}
 \phi_{st}(-y_i-m_i^{S^1}z+m_i^{(h_1,h_2)}(h_1)\tau-m_i^{(h_1,h_2)}(h_2),
  \tau,\sigma)
 \end{split} 
\end{equation} 
respectively. 

As for \eqref{eq:brvarA1} in Lemma \ref{lemm:brvarA} 
it is replaced by 
\begin{equation}\label{eq:brvarA3}
 \begin{split}
 & \quad \quad\brvar_{st}^A(X;z,\tau,\sigma)=(c\tau+d)^n 
 e^{-2\pi\img clz\sigma}\cdot \\
 &\quad \sum_{\Hatlam\in\HatLam}
 \frac{1}{|H_{\pi\circ\hatpi(\delta)}|}
 \int_{\HatF_{\Hatlam}}\sum_{(h_1,h_2)\in \delta}e(T'\HatF_{\Hatlam}) 
 e^{-2\pi\img m(h_1)dk}e^{-2n\pi\img(y+m^{S^1}z+m(h_2))ck}\cdot  \\
 &\prod_{i=1}^ne^{2\pi\img m_i^{(h_1,h_2)}(h_1)(c\tau+d)\sigma}
 \phi_{st}(-y_i-m_i^{S^1}z+m_i^{(h_1,h_2)}(h_1)\tau-m_i^{(h_1,h_2)}(h_2),
  \tau,(c\tau+d)\sigma).
\end{split} 
\end{equation}

Similarly \eqref{eq:hatvarA2} in {Lemma \ref{lemm:hatvarA} 
is replaced by
\begin{equation}\label{eq:hatvarA4}
 \begin{split}
 &\hatvar_{st}^A(X;z,\tau,\sigma) = (c\tau+d)^n
 e^{-2\pi\img clz\sigma}\sum_{\Hatlam\in\HatLam}
 \frac{1}{|H_{\pi\circ\hatpi(\delta)}|}
 \int_{\HatF_{\Hatlam}}\sum_{(h_1,h_2)\in \delta}e(T'\HatF_{\Hatlam}) 
  e^{-2\pi\img(y+m^{S^1}z)ck}\cdot \\
 &\prod_{i=1}^ne^{2\pi\img m_i^{(h_1,h_2)}(h_1)(c\tau+d)\sigma}
 \phi_{st}(-y_i-m_i^{S^1}z+m_i^{(h_1,h_2)}(h_1)\tau-m_i^{(h_1,h_2)}(h_2),
  \tau,(c\tau+d)\sigma).
\end{split} 
\end{equation}

By using \eqref{eq:stfxtptformula}, \eqref{eq:stbrfxtptformula}, 
\eqref{eq:brvarA3} and \eqref{eq:hatvarA4} proofs of stably 
almost complex versions of Theorems 
\ref{thm:brvar}, \ref{thm:hatvar} and \ref{thm:hatvartorsion} 
can be easily completed. 

Todd genus $T_0(X)$ is an invariant of stably almost complex orbifolds. 
It is defined as the index of a spin-c Dirac operator $D$ and written  
as the integral $T_0(X)=\int_XTd(T'X)$. We define the stabilized 
$T_y$-genus and stabilized $\hatT_y$-genus of $X$ by 
\[ \begin{split}
T_{y,st}(X)=&\ind(D\otimes\Lambda_y\tilde{T}^*X) \\
\hatT_{y,st}(X)=&\sum_{\hatlam\in\hatLam}(-y)^{f_{\hatlam}}
T_{y,st}(\hatX_{\hatlam}).
\end{split} \]
They are the degree zero terms in the $q$-expansions of 
$\zeta^{\frac{n}{2}}\var_{st}(X)$ and $\zeta^{\frac{n}{2}}
\hatvar_{st}(X)$.
If $N>1$ is an integer such that it is relatively prime to all 
$|H_x|$, then we define
\[ \brT_{y,st}(X)= \sum_{\hatlam\in\hatLam}(-y)^{\breve{f}_{\hatlam}}
T_{y,st}(\hatX_{\hatlam}). \]
Note that, if $X$ is an almost complex orbifold of dimension $2n$, 
then 
\[ T_{y,st}(X)=T_y(X)/(1+y)^n,\ \hatT_{y,st}(X)=\hatT_y(X)/(1+y)^n, \ 
 \brT_{y,st}(X)= \brT_y(X)/(1+y)^n. \]
We may define $T_y(X),\ \hatT_y(X) \ \text{and}\ \brT_y(X)$ for 
stably almost complex orbifolds using the above equalities. 
With these understandings the results in Section 7 still hold 
for stably almost complex orbifolds.

\providecommand{\bysame}{\leavevmode\hbox to3em{\hrulefill}\thinspace}

\end{document}